\documentclass[12pt]{amsart}
\usepackage{amssymb}
\usepackage[all]{xy}
\usepackage[toc,page,title,titletoc,header]{appendix}

\begin{document}
\theoremstyle{plain}
\newtheorem{thm}{Theorem}[section]
\newtheorem*{thm1}{Theorem 1}
\newtheorem*{thm1.1}{Theorem 1.1}
\newtheorem*{thmM}{Main Theorem}
\newtheorem*{thmA}{Theorem A}
\newtheorem*{thm2}{Theorem 2}
\newtheorem{lemma}[thm]{Lemma}
\newtheorem{lem}[thm]{Lemma}
\newtheorem{cor}[thm]{Corollary}
\newtheorem{pro}[thm]{Proposition}
\newtheorem{propose}[thm]{Proposition}
\newtheorem{variant}[thm]{Variant}
\theoremstyle{definition}
\newtheorem{notations}[thm]{Notations}
\newtheorem{rem}[thm]{Remark}
\newtheorem{rmk}[thm]{Remark}
\newtheorem{rmks}[thm]{Remarks}
\newtheorem{defi}[thm]{Definition}
\newtheorem{exe}[thm]{Example}
\newtheorem{claim}[thm]{Claim}
\newtheorem{ass}[thm]{Assumption}
\newtheorem{prodefi}[thm]{Proposition-Definition}
\newtheorem{que}[thm]{Question}
\newtheorem{con}[thm]{Conjecture}
\newtheorem*{con1.5*}{Conjecture 1.5*}

\newtheorem*{tvcon}{Tate-Voloch Conjecture}
\newtheorem*{dmmcon}{Dynamical Manin-Mumford Conjecture}
\newtheorem*{dmlcon}{Dynamical Mordell-Lang Conjecture}
\newtheorem*{condml}{Dynamical Mordell-Lang Conjecture}
\numberwithin{equation}{section}
\newcounter{elno}                
\def\points{\list
{\hss\llap{\upshape{(\roman{elno})}}}{\usecounter{elno}}}
\let\endpoints=\endlist

\newcommand{\GO}{{\rm GO}}
\newcommand{\Fan}{{(\F^{\an})}}
\newcommand{\Phian}{{(\Phi^{\an})}}
\newcommand{\lcm}{{\rm lcm}}
\newcommand{\Tor}{{\rm Tor}}
\newcommand{\perf}{{\rm perf}}
\newcommand{\ad}{{\rm ad}}
\newcommand{\Spa}{{\rm Spa}}
\newcommand{\Perf}{{\rm Perf}}
\newcommand{\alHom}{{\rm alHom}}
\newcommand{\SH}{\rm SH}
\newcommand{\Tan}{\rm Tan}
\newcommand{\res}{{\rm res}}
\newcommand{\Om}{\Omega}
\newcommand{\om}{\omega}
\newcommand{\la}{\lambda}
\newcommand{\mc}{\mathcal}
\newcommand{\mb}{\mathbb}
\newcommand{\surj}{\twoheadrightarrow}
\newcommand{\inj}{\hookrightarrow}
\newcommand{\zar}{{\rm zar}}
\newcommand{\Exc}{\rm Exc}
\newcommand{\an}{{\rm an}}
\newcommand{\red}{{\rm \mathbf{red}}}
\newcommand{\codim}{{\rm codim}}
\newcommand{\Supp}{{\rm Supp}}
\newcommand{\rank}{{\rm rank}}
\newcommand{\Ker}{{\rm Ker \ }}
\newcommand{\Pic}{{\rm Pic}}
\newcommand{\Div}{{\rm Div}}
\newcommand{\Hom}{{\rm Hom}}
\newcommand{\im}{{\rm im}}
\newcommand{\Spec}{{\rm Spec \,}}
\newcommand{\Nef}{{\rm Nef \,}}
\newcommand{\Frac}{{\rm Frac \,}}
\newcommand{\Sing}{{\rm Sing}}
\newcommand{\sing}{{\rm sing}}
\newcommand{\reg}{{\rm reg}}
\newcommand{\Char}{{\rm char}}
\newcommand{\Tr}{{\rm Tr}}
\newcommand{\ord}{{\rm ord}}
\newcommand{\id}{{\rm id}}
\newcommand{\NE}{{\rm NE}}
\newcommand{\Gal}{{\rm Gal}}
\newcommand{\Min}{{\rm Min \ }}
\newcommand{\Max}{{\rm Max \ }}
\newcommand{\Alb}{{\rm Alb}\,}
\newcommand{\GL}{{\rm GL}\,}        
\newcommand{\PGL}{{\rm PGL}\,}
\newcommand{\Bir}{{\rm Bir}}
\newcommand{\Aut}{{\rm Aut}}
\newcommand{\End}{{\rm End}}
\newcommand{\Per}{{\rm Per}\,}
\newcommand{\ie}{{\it i.e.\/},\ }
\newcommand{\niso}{\not\cong}
\newcommand{\nin}{\not\in}
\newcommand{\soplus}[1]{\stackrel{#1}{\oplus}}
\newcommand{\by}[1]{\stackrel{#1}{\rightarrow}}
\newcommand{\longby}[1]{\stackrel{#1}{\longrightarrow}}
\newcommand{\vlongby}[1]{\stackrel{#1}{\mbox{\large{$\longrightarrow$}}}}
\newcommand{\ldownarrow}{\mbox{\Large{\Large{$\downarrow$}}}}
\newcommand{\lsearrow}{\mbox{\Large{$\searrow$}}}
\renewcommand{\d}{\stackrel{\mbox{\scriptsize{$\bullet$}}}{}}
\newcommand{\dlog}{{\rm dlog}\,}    
\newcommand{\longto}{\longrightarrow}
\newcommand{\vlongto}{\mbox{{\Large{$\longto$}}}}
\newcommand{\limdir}[1]{{\displaystyle{\mathop{\rm lim}_{\buildrel\longrightarrow\over{#1}}}}\,}
\newcommand{\liminv}[1]{{\displaystyle{\mathop{\rm lim}_{\buildrel\longleftarrow\over{#1}}}}\,}
\newcommand{\norm}[1]{\mbox{$\parallel{#1}\parallel$}}
\newcommand{\boxtensor}{{\Box\kern-9.03pt\raise1.42pt\hbox{$\times$}}}
\newcommand{\into}{\hookrightarrow}
\newcommand{\image}{{\rm image}\,}
\newcommand{\Lie}{{\rm Lie}\,}      
\newcommand{\CM}{\rm CM}
\newcommand{\sext}{\mbox{${\mathcal E}xt\,$}}  
\newcommand{\shom}{\mbox{${\mathcal H}om\,$}}  
\newcommand{\coker}{{\rm coker}\,}  
\newcommand{\sm}{{\rm sm}}
\newcommand{\pgcd}{\text{pgcd}}
\newcommand{\trd}{\text{tr.d.}}
\newcommand{\tensor}{\otimes}
\renewcommand{\iff}{\mbox{ $\Longleftrightarrow$ }}
\newcommand{\supp}{{\rm supp}\,}
\newcommand{\ext}[1]{\stackrel{#1}{\wedge}}
\newcommand{\onto}{\mbox{$\,\>>>\hspace{-.5cm}\to\hspace{.15cm}$}}
\newcommand{\propsubset}
{\mbox{$\textstyle{
\subseteq_{\kern-5pt\raise-1pt\hbox{\mbox{\tiny{$/$}}}}}$}}
\newcommand{\sA}{{\mathcal A}}
\newcommand{\sB}{{\mathcal B}}
\newcommand{\sC}{{\mathcal C}}
\newcommand{\sD}{{\mathcal D}}
\newcommand{\sE}{{\mathcal E}}
\newcommand{\sF}{{\mathcal F}}
\newcommand{\sG}{{\mathcal G}}
\newcommand{\sH}{{\mathcal H}}
\newcommand{\sI}{{\mathcal I}}
\newcommand{\sJ}{{\mathcal J}}
\newcommand{\sK}{{\mathcal K}}
\newcommand{\sL}{{\mathcal L}}
\newcommand{\sM}{{\mathcal M}}
\newcommand{\sN}{{\mathcal N}}
\newcommand{\sO}{{\mathcal O}}
\newcommand{\sP}{{\mathcal P}}
\newcommand{\sQ}{{\mathcal Q}}
\newcommand{\sR}{{\mathcal R}}
\newcommand{\sS}{{\mathcal S}}
\newcommand{\sT}{{\mathcal T}}
\newcommand{\sU}{{\mathcal U}}
\newcommand{\sV}{{\mathcal V}}
\newcommand{\sW}{{\mathcal W}}
\newcommand{\sX}{{\mathcal X}}
\newcommand{\sY}{{\mathcal Y}}
\newcommand{\sZ}{{\mathcal Z}}
\newcommand{\A}{{\mathbb A}}
\newcommand{\B}{{\mathbb B}}
\newcommand{\C}{{\mathbb C}}
\newcommand{\D}{{\mathbb D}}
\newcommand{\E}{{\mathbb E}}
\newcommand{\F}{{\mathbb F}}
\newcommand{\G}{{\mathbb G}}
\newcommand{\HH}{{\mathbb H}}
\newcommand{\I}{{\mathbb I}}
\newcommand{\J}{{\mathbb J}}
\newcommand{\M}{{\mathbb M}}
\newcommand{\N}{{\mathbb N}}
\renewcommand{\P}{{\mathbb P}}
\newcommand{\Q}{{\mathbb Q}}
\newcommand{\R}{{\mathbb R}}
\newcommand{\T}{{\mathbb T}}
\newcommand{\U}{{\mathbb U}}
\newcommand{\V}{{\mathbb V}}
\newcommand{\W}{{\mathbb W}}
\newcommand{\X}{{\mathbb X}}
\newcommand{\Y}{{\mathbb Y}}
\newcommand{\Z}{{\mathbb Z}}
\newcommand{\fX}{\mathfrak{X}}
\newcommand{\wF}{\widetilde{F}}
\newcommand{\wG}{\widetilde{G}}
\newcommand{\fix}{\mathrm{Fix}}

\title[]{Algebraic dynamics of the lifts of Frobenius}

\author{Junyi Xie}

\address{IRMAR
campus de Beaulieu, 
b\^atiments 22 et 23
263 avenue du G\'en\'eral Leclerc, CS 74205
35042  Rennes c\'edex}


\thanks{The author is supported by the labex CIMI}

\email{junyi.xie@univ-rennes1.fr}
\date{\today}
\bibliographystyle{plain}

\maketitle

\maketitle

\begin{abstract}
We study the algbraic dynamics for endomorphisms of projective spaces with coefficients in a p-adic field
whose reduction in positive characteritic is the Frobenius. In particular, we prove a version of the dynamical Manin-Mumford conjecture and the dynamical Mordell-Lang conjecture for the coherent backward orbits for such endomorphisms. We also give a new proof of a dynamical version of the Tate-Voloch conjecture in this case.
Our method is based on the theory of perfectoid spaces introduced by P. Scholze. In the appendix, we prove that under some technical condition on the field of definition, a dynamical system for a polarized lift of Frobenius on a projective variety can be embedding into a dynamical system for some endomorphism of a projective space.
\end{abstract}

\tableofcontents

\newpage

\section{Introduction}
In this paper, we write $\C_p$ for the completion of the algebraically closure of $\Q_p$ with the induced norm.
Denote by $\C_p^{\circ}$ its valuation ring and $\C_p^{\circ\circ}$ the maximal ideal of $\C_p^{\circ}.$
Let $F:\P^N_{\C_p}\rightarrow \P^N_{\C_p}$ be an endomorphism taking form $$F:[x_0:\cdots:x_N]\mapsto [x_0^q+p'P_0(x_0,\cdots,x_N):\cdots: x_N^q+p'P_N(x_0,\cdots,x_N)]$$ where $q$ is a power of $p$, $p'\in \C_p^{\circ\circ},$ and 
$P_0,\cdots,P_N$ are homogeneous polynomials of degree $q$ in $\C_p^{\circ}[x_0,\cdots,x_N].$ 
We say that $F$ is \emph{a lift of Frobenius on $\P^N_{\C_p}$}.
%

In this paper we present a new argument to study the algebraic dynamics for such maps, which is based on the theory of perfectoid spaces introduced by  Scholze. In particular, we study some dynamical analogues of diophantine geometry for such maps.

\bigskip
\subsubsection*{Dynamical Manin-Mumford conjecture}
At first, we recall the dynamical Manin-Mumford conjecture proposed by Zhang \cite{Zhang1995}.

\begin{dmmcon}Let $F:X_{\C}\rightarrow X_{\C}$ be an endomorphism of a quasi-projective variety defined over $\C$.
Let $V$ be a subvariety of $X$.
If the Zariski closure of the set of preperiodic\footnote{A preperiodic point $x$ is a point satisfying $F^m(x)=F^n(x)$ for some $m>n\geq 0.$} (resp. periodic\footnote{A periodic point $x$ is a point satisfying $F^n(x)=x$ for some $n> 0$.} )
points of $F$ contained in $V$ is
Zariski dense in $V$, then $V$ itself is preperiodic (resp. periodic).
\end{dmmcon}

This conjecture is a dynamical analogue of
the Manin-Mumford conjecture on subvarieties of abelian varieties.
More precisely, let $V$ be an irreducible subvariety inside an abelian variety $A$ over $\C$ such that the intersection of the set of torsion points of $A$ and $V$ is Zariski dense in $V$.
Then the Manin-Mumford conjecture asserts that there exists an abelian subvariety $V_0$ of $A$ and a torsion point $a\in A(\C)$ such that $V=V_0+a$.

The Manin-Mumford conjecture was first proved by Raynaud \cite{Raynaud1983,Raynaud1983a}.
Various versions of this conjecture were proved
by Ullmo \cite{Ullmo1998}, Zhang \cite{Zhang1998}, Buium \cite{Buium1996}, Hrushovski \cite{Hrushovski2001} and Pink-Roessler\cite{Pink2002}.
Observe that the dynamical Manin-Mumford conjecture for the map $x\mapsto 2x$ on $A$ implies the classical Manin-Mumford conjecture.

\medskip

The dynamical Manin-Mumford conjecture does not hold in full generality, as we have some counterexamples \cite{Ghioca2011,Pazuki2010,Pazuki2013}.  In particular, Pazuki \cite{Pazuki2013} shows that counterexamples can come from a lift of Frobenius crossed with a lift of its Verschiebung.
This motivated the proposal of several modified versions of it \cite{Ghioca2011,Yuan}.

\medskip

However,
this conjecture is now known to hold in some special cases \cite{Baker2005,Fakhruddin,Medvdev,Ghioca2010,Ghioca2011,Dujardin,DragosGhioca, DragosGhiocaa}. It seems that the dynamical Manin-Mumford conjecture may be true except a few families of counterexamples.

\medskip
In this paper, we prove the dynamical Manin-Mumford conjecture for periodic points of lifts of Frobenius on $\P^N$.

\begin{thm}\label{thmmmdlfper}Let $F:\P^N_{\C_p}\rightarrow \P^N_{\C_p}$ be
a lift of Frobenius on $\P^N_{\C_p}$. Denote by $\Per$ the set of periodic closed points in $\P^N_{\C_p}$. Let $V$ be any irreducible subvariety of $\P^N_{\C_p}$ such that $V\cap\Per$ is Zariski dense in $V$. Then $V$ is periodic i.e. there exists ${\ell}\geq 1$ such that $F^{\ell}(V)=V.  $
\end{thm}

We note that in \cite{Medvdev}, Medvedev and Scanlon  have proved Theorem \ref{thmmmdlfper} in the case $$F:[x_0:\cdots:x_N]\mapsto [x_0^q+pP(x_0,x_N):\cdots: x_{N-1}^q+pP(x_{N-1},x_N):x_N^q]$$ where $q$ is a power of $p$ and $P\in \Z_p[x,y]$ is a homogenous polynomial of degree $q$. In \cite{Pazuki2013}, Pazuki have studied the lifts of Frobenius on abelian varieties.

\smallskip

We should mention that, recently Scanlon got a new proof of this theorem without using pefectoid spaces. Since this proof is unpublished and it is completely different from ours, we will discuss it briefly in Section \ref{sectionperiodic} of this paper.

\medskip

\subsubsection*{Dynamical Tate-Voloch Conjecture}
Let $V$ be an irreducible subvariety of $\P^N_{\C_p}$. There are homogenous polynomials $H_i\in \C_p[x_0,\dots,x_N]$, $i=1,\dots, m$ satisfying $\|H_i\|=1$ which define $V$. For any point $y\in \P^N_{\C_p}(\C_p)$, we may write $y=[y_0:\dots:y_N]$, $\max\{|y_i|\}_{0\leq i\leq N}=1.$ Then we denote by $d(y,V):=\max\{|H_i(y_0,\dots,y_N)|\}_{1\leq i\leq m}.$
Observe that $d(y,V)$ does not depend on the choice of $\{H_i\}_{1\leq i\leq m}$ and the coordinates $[y_0:\dots:y_N]$ of $y.$ It can be viewed as the distance between $y$ and $V$. Moreover for any quasi-projective variety $X$ and subvariety $V$ of $X$, by choosing an embedding $X\hookrightarrow  \P^N_{\C_p}$, $d(\bullet,\overline{V})$ defines a distance between $V$ and a point in $X.$

\medskip

In \cite{Tate1996}, Tate and Voloch made the following conjecture.
\begin{tvcon}
Let $A$ be a semiabelian variety over $\C_p$ and $V$ a subvariety of $A$. Then there exists $c>0$ such that for any torsion point $x\in A$, we have either $x\in V$ or $d(x,V)>c.$
\end{tvcon}
This conjecture is proved by Scanlon in \cite{Scanlon1999} when $A$ is defined over a finite extension of $\Q_p.$

In \cite{Buium1996a}, Buium proved a dynamical version of this conjecture for periodic points of lifts of Frobenius on any algebraic variety. Here we state it only for the lifts of Frobenius on $\P^N_{\C_p}.$
\begin{thm}\label{thmdvtperlf}Let $F:\P^N_{\C_p}\rightarrow \P^N_{\C_p}$ be
a lift of Frobenius on $\P^N_{\C_p}$. Let $V$ be any irreducible subvariety of $\P^N_{\C_p}$. Then there exists $\delta>0$ such that for any point $x\in \Per$, either $d(x,V)>\delta$ or $x\in V.$
\end{thm}

In this paper, we give a new proof of this theorem by using the theory of perfectoid spaces.

\medskip

\subsubsection*{Dynamical Mordell-Lang conjecture}
The Mordell-Lang conjecture on subvarieties of semiabelian varieties (now a theorem of Faltings \cite{Faltings1994} and Vojta \cite{Vojta1996}) says that if $V$ is a subvariety of a semiabelian variety $G$ defined over $\C$ and $\Gamma$ is a finitely generated subgroup  of $G(\C)$, then $V(\C)\bigcap \Gamma$ is a union of at most finitely many translates of
subgroups of $\Gamma$.

Inspired by this, Ghioca and Tucker proposed the following dynamical analogue of the Mordell-Lang conjecture.
\begin{dmlcon}[\cite{Ghioca2009}]\label{dml}Let $X$ be a quasi-pro\-jective variety defined over $\mathbb{C}$, let
$f: X \rightarrow X$ be an endomorphism, and $V$ be any subvariety of $X$. For any point $x\in X(\mathbb{C})$  the set $\{n\in \mathbb{N}|\,\,f^n(x)\in V(\mathbb{C})\}$ is a union of at most finitely many arithmetic progressions\footnote{An arithmetic progression is a set of the form $\{an+b|\,\,n\in \mathbb{N}\}$ with $a,b\in \mathbb{N}$. In particular, when $a=0$, it contains only one point}.
\end{dmlcon}

Observe that the dynamical Mordell-Lang conjecture implies the classical Mordell-Lang conjecture in the case $\Gamma\simeq (\mathbb{Z},+)$.

\medskip

The dynamical Mordell-Lang conjecture is proved in many cases. For example, Bell, Ghioca and Tucker proved this conjecture for \'etale maps \cite{Bell2010}, and the author proved it for the endomorphisms on $\A^2_{\bar{\Q}}$ \cite{Xiec}. We refer to the book \cite{Bell2016} for a good survey of this conjecture.

We note that the dynamical Mordell-Lang conjecture is not a full generalization of the Mordell-Lang conjecture. In particular, it considers only the forward orbit but not the backward orbit.  In an informal seminar, S-W Zhang asked the following question to me.
\begin{que}\label{quezhang}Let $X$ be a quasi-projective variety over $\C$ and $F:X\to X$ be a finite endomorphism. Let $x$ be a point in $X(\C)$. Denote by $O^-(x):=\cup^{\infty}_{i=0}F^{-i}(x)$ the backward orbit of $x$. Let $V$ be a positively dimensional irreducible subvariety of $X$. If $V\cap O^-(x)$ is Zariski dense in $V$, what can we say about $V$?
\end{que}

We note that if $V$ is preperiodic, then $V\cap O^-(x)$ is Zariski dense in $V$. As the dynamical Manin-Mumford conjecture, the converse is not true.  Indeed, we have the following example. Let $X=\A^1_{\C}\times \A^1_{\C}$ and $f:X\to X$ be the endomorphism defined by $(x,y)\mapsto (x^4,y^6).$ Let $V$ be the diagonal and $x=(1,1)$. Then $V\cap O^-(x)$ is Zariski dense in $V$, but $V$ is not preperiodic. 
We have counterexamples even when $F$ is a polarized\footnote{An endomorphism $F:X\to X$ on a projective variety is said to be  polarized if there exists an ample line bundle $L$ on $X$ satisfying $F^*L=L^{\otimes d}, d\geq 2$.} endomorphism. The following example is given by Ghioca, which is similar to \cite[Theorem 1.2]{Ghioca2011}.
\begin{exe}Let $E$ be the elliptic curve over $\C$ defined by the lattice $\Z[i]\subseteq \C.$
Let $F_1$ be the endomorphism on $E$ defined by the multiplication by $10$ and $F_2$ be the endomorphism on $E$ defined by the multiplication by $6+8i$.
Set $X:=E\times E$, $F:=(F_1,F_2)$ on $X$. Since $|10|=|6+8i|$, $F$ is a polarized endomorphism on $X.$ 
Let $V$ be the diagonal in $X$ and $x$ be the origin. We may check that $V\cap O^-(x)$ is Zariski dense in $V$, but $V$ is not preperiodic. 
\end{exe}

As a special case of Question \ref{quezhang}, we propose the following conjecture.
\begin{con}\label{condmlrev}Let $X$ be a quasi-projective variety over $\C$ and $F:X\to X$ be a finite endomorphism. Let $\{b_i\}_{i\geq 0}$ be a sequence of points in $X(\C)$ satisfying $f(b_i)=b_{i-1}$ for all $i\geq 1$. Let $V$ be a positively dimensional irreducible subvariety of $X$. If the $\{b_i\}_{i\geq 0}\cap V$ is Zariski dense in $V$, then $V$ is periodic under $F$.
\end{con}
\begin{rem}This conjecture can be viewed as the dynamical Mordell-Lang conjecture for the coherent backward orbits. In fact, it is easy to see that Conjecture \ref{condmlrev} is equivalent to the following:

\begin{con1.5*}
Let $X$ be a quasi-projective variety over $\C$ and $F:X\to X$ be a finite endomorphism. Let $\{b_i\}_{i\geq 0}$ be a sequence of points in $X(\C)$ satisfying $f(b_i)=b_{i-1}$ for all $i\geq 1$. Let $V$ be a subvariety of $X$. Then the set $\{n\geq 0|\,\, b_n\in V\}$ is a union of at most finitely many arithmetic progressions.
\end{con1.5*}

\proof[Conjeture \ref{condmlrev}$\Rightarrow$ Conjecture 1.5*]
Set $W:=\cap_{n\geq 0}\overline{\{b_i|\,\, b_i\in V, i\geq n\}}.$ Then there exists $N\geq 0$ such that $W=\overline{\{b_i|\,\, b_i\in V, i\geq N\}}$. 
We note that $\{n\geq 0|\,\, b_n\in V\}\setminus \{n\geq 0|\,\, b_n\in W\} \subseteq \{0,\dots,N\}$ is finite.
After replacing $b_0$ by $b_N$, we may assume that $N=0$. 
If $W$ is empty, then $\{n\geq 0|\,\, b_n\in V\}=\{n\geq 0|\,\, b_n\in W\}=\emptyset$.
If $W$ is not empty, then every irreducible component of $W$ has positive dimension and $\{b_i\}_{i\geq N}\cap W$ is Zariski dense in $W$. 
Conjeture \ref{condmlrev} implies that there exists $r\geq 1$ such that $F^r(W)=W$. If for some index $i\in \{0,\dots,r-1\}$, there exists $s\geq 0$ such that $b_{i+sr}\not\in V$, then $b_{i+nr}\not\in V$ for all $n\geq s$. Denote by $T_i, i=0,\dots,r-1$ the set of $j\geq 0$ satisfying $b_j\in V$ and $j=i \mod r.$ Then $T_i$ is either finite or equal to $\{i+rn|\,\, n\in \N\}.$ It follows that $$\{n\geq 0|\,\, b_n\in V\}=\{n\geq 0|\,\, b_n\in W\}=\cup_{i=0}^{r-1}T_i$$ is a union of at most finitely many arithmetic progressions.

\proof[Conjeture \ref{condmlrev}$\Rightarrow$ Conjeture 1.5*]
Assume that $V$ is a positively dimensional irreducible subvariety of $X$ such that $\{b_i\}_{i\geq 0}\cap V$ is Zariski dense in $V$.
Then $\{n\geq 0|\,\, b_n\in V\}$ is infinite.
Conjeture \ref{condmlrev} shows that $\{n\geq 0|\,\, b_n\in V\}$ takes form $\{n\geq 0|\,\, b_n\in V\}=F\cup(\cup_{j=1}^sT_j)$ where $F$ is finite and $T_j, j=1,\dots,s$ are infinite arithmetic progressions. There exists $j\in \{1,\dots,s\}$ such that $\{b_i|\,\, i\in T_j\}$ is Zariski dense in $V.$
Write $T_j=a+r\N$ where $a\geq 0, r\geq 1$. Since $F(\{b_i|\,\, i\in T_j\})\setminus \{b_i|\,\, i\in T_j\}=\{a\}$, we have $F^r(V)=V.$
\end{rem}

 In this paper, we prove Conjecture \ref{condmlrev} for the lifts of Frobenius of $\P^N_{\C_p}.$
\begin{thm}\label{thmrevedml}Let $F:\P^N_{\C_p}\rightarrow \P^N_{\C_p}$ be
a lift of Frobenius on $\P^N_{\C_p}$. Let $\{b_i\}_{i\geq 0}$ be a sequence of points in $\P^N_{\C_p}(\C_p)$ satisfying $f(b_i)=b_{i-1}$ for all $i\geq 1$. Let $V$ be a positively dimensional irreducible subvariety of $\P^N_{\C_p}$. If the $\{b_i\}_{i\geq 0}\cap V$ is Zariski dense in $V$, then $V$ is periodic under $F$.
\end{thm}

In fact, we prove a stronger statement.
\begin{thm}\label{thmrevedmllimit}Let $F:\P^N_{\C_p}\rightarrow \P^N_{\C_p}$ be
a lift of Frobenius on $\P^N_{\C_p}$. Let $\{b_i\}_{i\geq 0}$ be a sequence of points in $\P^N_{\C_p}(\C_p)$ satisfying $f(b_i)=b_{i-1}$ for all $i\geq 1$. Let $V$ be a subvariety of $\P^N_{\C_p}$. If there exists a subsequence $\{b_{n_i}\}_{i\geq 0}$ such that $|d(b_{n_i},V)|\to 0$ when $n\to\infty$, then $b_{n_i}\in V$ for $i$ large enough and there exists $r\geq 0$, such that $\{b_i\}_{i\geq 0}\subseteq \cup_{i=0}^rF^i(V).$
\end{thm}
It implies the following Tate-Voloch type statement.
\begin{cor}\label{cordmlrevtv}Let $F:\P^N_{\C_p}\rightarrow \P^N_{\C_p}$ be
a lift of Frobenius on $\P^N_{\C_p}$. Let $\{b_i\}_{i\geq 0}$ be a sequence of points in $\P^N_{\C_p}(\C_p)$ satisfying $f(b_i)=b_{i-1}$ for all $i\geq 1$. Let $V$ be a subvariety of $\P^N_{\C_p}$. Then there exists $c>0$ such that for all $i\geq 0$, either $b_i\in V$ or $d(b_i,V)>c.$
\end{cor}

\bigskip

\subsubsection*{Overview of the proofs}

Let us now see in more detail how our arguments work. 

Denote by $K:=\C_p$ and $K^{\flat}:=\widehat{\overline{\F_p((t))}}$ the completion of the algebraic closure of $\F_p$.
We denote by $K^{\circ}$ (resp. $K^{\flat\circ}$) the valuation ring of $K$ (resp. $K^{\flat}$) and denote by $K^{\circ\circ}$ (resp. $K^{\flat\circ\circ}$) the maximal ideal of  $K^{\circ}$ (resp. $K^{\flat\circ}$).
Denote by $k:=\overline{\F_p}$. We have $k=K^{\circ}/K^{\circ\circ}=K^{\flat\circ}/K^{\flat\circ\circ}.$ Moreover, we have an embedding $k\hookrightarrow K^{\flat}.$

Let $F:\P^N_{K}\rightarrow \P^N_{K}$ be an endomorphism taking form $$F:[x_0:\cdots:x_N]\mapsto [x_0^q+p'P_0(x_0,\cdots,x_n):\cdots: x_N^q+p'P_N(x_0,\cdots,x_N)]$$ where $p'\in K^{\circ\circ}$, $q$ is a power of $p$, and $P_0,\cdots,P_N$ are homogeneous polynomials of degree $q$ in $K^{\circ}[x_0,\cdots,x_N].$

 We associate to $\P^N_{K}$ ( resp. $\P^N_{K^{\flat}}$) a non archimedean analytic space $\P^{N,\ad}_{K}$ (resp. $\P^{N,\ad}_{K^{\flat}}$) with a natural embedding $\P^N_{K}(K)\subseteq \P^{N,\ad}_{K}$ (resp. $\P^N_{K^{\flat}}(K^{\flat})\subseteq\P^{N,\ad}_{K^{\flat}}$). The endomorphism $F$ extends to an endomorphism $F^{\ad}$ on $\P^{N,\ad}_{K}.$

 Denote by $\lim\limits_{\overleftarrow{F^{\ad}}}\P^{N,\ad}_{K}$ the inverse limit $\P^{N,\ad}_{K}$'s where transition maps are $F^{\ad}.$
 Then we may construct a perfectoid space $\P^{N,\perf}_{K}$ with an endomorphism $F^{\perf}$ for which the topological dynamical system $(\P^{N,\perf}_{K}, F^{\perf})$ is isomorphic to $(\lim\limits_{\overleftarrow{F^{\ad}}}\P^{N,\ad}_{K},T)$ where $T:(x_0,x_1,\cdots)\to (F^{\ad}(x_0),x_0,\cdots)$ is the shift map on $\lim\limits_{\overleftarrow{F^{\ad}}}\P^{N,\ad}_{K}$. Moreover, we have a natural morphism $\pi:\P^{N,\perf}_{K}\to \P^{N,\ad}_{K}$ defined by the projection to the first coordinate.
 This construction has been stated by Scholze \cite[Section 7]{Scholze2014}.

 Similarly, we construct a perfectoid space $\P^{N,\perf}_{K^{\flat}}$ which is isomorphic to the inverse limit $\lim\limits_{\overleftarrow{\Phi^s}}\P^{N,\ad}_{K^{\flat}}$ where $\Phi$ is the Frobenius endomorphism on $\P^{N,\ad}_{K^{\flat}}$, and $q=p^s$. Denote by $\pi^{\flat}:\P^{N,\perf}_{K^{\flat}}\to \P^{N,\ad}_{K^{\flat}}$ the morphism defined by the projection to the first coordinate. Since $\Phi$ is a homeomorphism on the underlying topological space, $\pi^{\flat}$ induces an isomorphism from the topological dynamical system $(\P^{N,\perf}_{K^{\flat}},\Phi^{s,\perf})$ to $(\P^{N,\ad}_{K^{\flat}},\Phi^{s,\ad})$, where $\Phi^{\perf}$ is the Frobenius on $\P^{N,\perf}_{K^{\flat}}.$

By the theory of perfectoid spaces, there is a natural homeomorphism of topological space $\rho:\P^{N,\perf}_{K}\to \P^{N,\perf}_{K^{\flat}}$ satisfying $\Phi^{s,\perf}\circ\rho=\rho\circ F^{\perf}.$

\medskip

As an example, we explain the proof of Theorem \ref{thmmmdlfper}.
Let $V$ be any subvariety of $\P^N_{\C_p}$ such that $V\cap\Per$ is Zariski dense in $V$.

It is easy to see that the map $\pi\circ\rho^{-1}\circ(\pi^{\flat})^{-1}$ induces a bijection from the set $\Per^{\flat}$  of periodic points of $\Phi^s$ in $\P^N_{K^{\flat}}(K^{\flat})$ to the set $\Per$  of periodic points of $F$ in $\P^N_K(K).$ We note that the set of periodic points of $\Phi^s$ in $\P^N_{K^{\flat}}(K^{\flat})$ is exactly the set of points defined over $k$ i.e. the image of $\eta:\P^N_k(k)\hookrightarrow \P^N_{K^{\flat}}(K^{\flat}).$ We have a reduction map $\red: \P^N_K(K)\to \P^N_k(k)$. The map $\eta\circ\red:\Per\to \Per^{\flat}$ is bijective. Moreover, we have that $(\eta\circ\red)\circ(\pi\circ\rho^{-1}\circ(\pi^{\flat})^{-1})$ is identity on $\Per^{\flat}$.

Denote by $S^{\flat}$ the Zariski closure of $\eta\circ\red(V\cap \Per).$ Since $S^{\flat}$ is defined over $k$, it is periodic under $\Phi^s$. The main ingredient of our proof is to show that $S^{\flat}$ is a subset of $\pi^{\flat}(\rho(\pi^{-1}(V))).$ If $\pi^{\flat}(\rho(\pi^{-1}(V)))$ is algebraic, this is obvious. But a priori, $\pi^{\flat}(\rho(\pi^{-1}(V)))$ is not algebraic, since the map $\rho$ is very transcendental. Our strategy is to approximate $\pi^{\flat}(\rho(\pi^{-1}(V)))$ by algebraic subvarieties of $\P^N_{K^{\flat}}$. For simplicity, assume that $V$ is an hypersurface of $\P^N_{K}$.  Applying the approximation lemma \cite[Corollary 6.7]{Scholze2012} of Scholze,  for any $\epsilon>0$, there exists an algebraic hypersurface $H_{\epsilon}$ of $\P^N_{K^{\flat}}$ which is $\epsilon$-close to $\pi^{\flat}(\rho(\pi^{-1}(V)))$. Then $\eta\circ\red(V\cap \Per)$ is $\epsilon$-close to $H_{\epsilon}$. Since $S^{\flat}$ the Zariski closure of $\eta\circ\red(V\cap \Per)$ in $\P^N_{K^{\flat}}$, we can show that it is $\epsilon$-close to $H_{\epsilon}$. Then we can show that $S^{\flat}$ is contained in $\pi^{\flat}(\rho(\pi^{-1}(V)))$ by letting $\epsilon$ tends to $0$.

Then we have $S:=\pi(\rho^{-1}((\pi^{\flat})^{-1}(S^{\flat})))\subseteq V$. Since $S$ is periodic and Zariski dense in $V$, then $V$ is periodic.

%
%
%

\bigskip


%
In this paper, we mainly consider the lifts of Frobenius on $\P^N_{\C_p}$ for the simplicity, since we think that the aim of this paper is to present this new method in dynamics. 
We suspect that our method can be applied to the more general case where $F$ is a lift of Frobenius on any projective variety over $\C_p$. On the other hand, it is often that a lift of Frobenius on a projective variety $X$ can be extended to some lift of Frobenius on $\P^N_{\C_p}$ for some embedding $\tau:X\hookrightarrow \P^N_{\C_p}.$ In the appendix, we prove the existence of such embedding for \emph{polarized} lifts of Frobenius on some projective varieties under some technical condition on the field of definition.
Once this happens, many questions can be reduced to the special case where $X=\P^N_{\C_p}.$

\bigskip

\subsubsection*{The plan of the paper}
The paper is organized as follows. In Section \ref{sectionpre}, we gather a number of results on the perfectoid spaces in Scholze's papers \cite{Scholze2012,Scholze2014}.
In Section \ref{sectioninvlim}, we construct the inverse limit and make it to be a perfectoid space with an automorphism. We also construct its tilt and give the isomorphism between these two topological dynamical systems. In Section \ref{sectionperiodic}, we study the periodic points of $F$. In particular, we prove  Theorem \ref{thmmmdlfper} and Theorem \ref{thmdvtperlf}. In Section \ref{sectioninverseorbit}, we study the coherent backward orbits of a point. In particular, we prove Theorem \ref{thmrevedml}, Theorem \ref{thmrevedmllimit} and Corollary \ref{cordmlrevtv}.  In the apprendix,  we study the polaried lift of Frobenius on projective varieties over $\C_p$.

%

\section*{Acknowledgement}
I would like to thank Charles Favre, Serge Cantat, St\'ephane Lamy, Jean Gillibert, Dragos Ghioca, Thomas Tucker  and Peter Scholze for useful discussions.
I thank Fabien Mehdi Pazuki and Umberto Zannier for their comments on the first version of this paper. I thank Thomas Scanlon, who told me about his new proof of Theorem \ref{thmmmdlfper} and let me know the Tate-Voloch conjecture. We thank the referees for numerous insightful remarks.
I especially thank Shou-Wu Zhang, who introduced the theory of perfectoid spaces to me and posed me Question \ref{quezhang}.

\newpage

\section{Preliminary: perfectoid spaces}\label{sectionpre}
In this section, we introduce some necessary background in perfectoid spaces. All the results in this section can be found in Scholze's papers \cite{Scholze2012,Scholze2014}.
The perfectoid spaces are some nonarchimedean analytic spaces. Following the technique of Scholze in \cite{Scholze2012}, we work with Huber's adic spaces \cite{Huber1993,Huber1994,Huber1996}.
\subsection{Adic spaces}
In this section, we denote by $k$ a complete nonarchimedean field i.e. a complete topological field whose topology is induced by a nontrivial norm $|\cdot|:k\to [0,\infty).$ Denote by $R$ a topological $k$-algebra. Moreover we suppose that $R$ is a Tate $k$-algebra i.e.  there exists a subring $R_0\subseteq R$, such that $aR_0$, $a\in k^{\times}$, forms a basis of open neighborhoods of $0$.

A subset $M\subseteq R$ is call bounded if $M\subseteq aR_0$, for some $a\in k^{\times}.$ An element $x\in R$ is called power-bounded if $\{x^n|\,\,n\geq 0\}\subseteq R$ is bounded. Let $R^{\circ}\subseteq R$ be the subring of power-bounded elements.

\begin{defi}[\cite{Scholze2012}]An affinoid $k$-algebra is a pair $(R,R^+)$, where $R$ is a Tate $k$-algebra and $R^+$ is an open and integrally closed subring of $R^{\circ}.$
\end{defi}

A valuation on $R$ is a map $|\cdot|: R\to \Gamma\cup \{0\}$, where $\Gamma$ is a totally ordered abelian group, such that, $|0|=0$, $|1|=1$, $|xy|=|x||y|$ and $|x+y|\leq \max\{|x|,|y|\}.$ We say that $|\cdot|$ is continuous, if for all $\gamma\in \Gamma$, the subset $\{x\in R|\,\, |x|<\gamma\}\subseteq R$ is open.

To a pair $(R,R^+)$, Huber associates a space $\Spa(R,R^+)$ of equivalence classes of continuous valuations $|\cdot|$ on $R$ such that $|R^{+}|\leq 1,$ and call it an \emph{affinoid space}.

For a point $x\in \Spa(R,R^+)$, we denote by $f\to |f(x)|$ the associated valuation.
It is a fact \cite[Proposition 2.12. (iii)]{Scholze2012} that $$R^+=\{f\in R|\,\,|f(x)|\leq 1 \text{ for all } x\in \Spa(R,R^+)\}.$$
We equip $\Spa(R,R^+)$ with the topology generated by rational subsets:
$$
U(f_1,\cdots,f_n;g)=\{x\in \Spa(R,R^+)|\,\, |f_i(x)| \leq |g(x)| \}\subseteq \Spa(R,R^+)
$$
where $f_1,\cdots,f_n\in R$ generate $R$ as an idea and $g\in R$.

The completion $(\widehat{R},\widehat{R}^+)$ of an affinoid algbra $(R,R^+)$ is also an affinoid algebra. Then we recall \cite[Proposition 3.9]{Huber1993}.
\begin{pro}\label{procompletionaffi}We have $\Spa(\widehat{R},\widehat{R}^+)\simeq \Spa(R,R^+)$, identifying rational subsets.
\end{pro}

We say a point $x\in \Spa(R,R^+)$ is a $k$-point, if the valuation $x$ is induced by a morphism from $R$ to $k$ i.e. there exists a morphism $\phi: R\to k$ such that for any $f\in R$, $|f(x)|=\phi(f)|$.

\medskip

Roughly speaking,  \emph{adic spaces} over $K$ are the objects obtained by gluing affinoid spaces. 
The morphisms betweens the adic spaces are the morphisms glued by the morphisms between affinoid spaces.  Because in this paper we only consider some very concrete adic spaces,  we just give a very brief definition of the adic spaces.  One may find a detailed definition in \cite{Huber1994}.

\medskip

On an affinoid space $X=\Spa(R,R^+)$, one may define a pre-sheaves $O_X$ and $O_X^{+}$ on $X$. Since we do not use these pre-sheaves in this paper, we omit its definition.  
We do \emph{not} know whether $O_X$ is a sheaf in general. We note that once $O_X$ is a sheaf, $O_X^{+}$ is a sheaf also.  However, if $(R,R^+)$ is of \emph{topological finite type}\footnote{An affinoid $k$-algebra $(R,R^+)$ is said to be  of topological finite type if $R$ is a quotient of $k\{T_1,\dots,T_n\}$ for some $n$, and $R^+=R^{\circ}$.} then $O_X$ is a sheaf.  Assume that $O_X$ is a sheaf on $X$. For any $x\in X$, the valuation $f\mapsto |f(x)|$ extends to the stalk $O_{X,x}$, and we have $O_{X,x}^+=\{f\in O_{X,x}|\,\, |f(x)|\leq 1\}.$
The affinoid spaces $X$ defines a triple $(X, O_X, |\cdot (x)| \,|\,\, x\in X).$ 

An adic space over $k$ is a triple $(Y, O_Y, |\cdot (x)| \,|\,\, x\in Y)$, consisting of  a locally ringed topological space $(Y, O_Y)$ where $O_Y$ is a sheaf of complete topological $k$-algebras, and a continuous valuation $|\cdot (x)|$ on $O_{X,x}$ for every $x\in X$, which is locally on $Y$ an affinoid adic space.

\medskip

Let $X$ be an affinoid space. We say a point $x\in X$ is a $k$-point if it is a $k$-point in any (and thus all) affinoid neighbourhood of $X.$

%
%
%
%
%
%

\subsection{Perfectoid fields}
Denote by $K$ a complete nonarchimedean field of residue characteristic $p>0$ with norm $|\cdot|:K\to \R_{\geq 0}$.
Denote by $K^{\circ}:=\{x\in K| \,\,|x|\leq 1\}$ its valuation ring.

\begin{defi} We say $K$ is a \emph{perfectoid field} if $|K|\subseteq \R_{\geq 0}$ is dense in $\R_{\geq 0}$ and the Frobenius map $\Phi: K^{\circ}/p\to K^{\circ}/p$ is surjective.
\end{defi}

Observe that $\C_p$ and $\widehat{\overline{\F_p((t))}}$ are perfectoid fields.
Set $\Q_p(p^{1/p^{\infty}}):=\cup_{i\geq 0}\Q_p(p^{1/p^{i}})$ and $\F_p((t))(t^{1/p^{\infty}}):=\cup_{i\geq 0}\F_p((t))(t^{1/p^{i}})$. Then their completions $\widehat{\Q_p(p^{1/p^{\infty}})}$ and $\widehat{\F_p((t))(t^{1/p^{\infty}})}$ are perfectoid fields. Note that $\Q_p$ is not a perfectoid field, since $|\Q_p|={0}\cup \{p^i|\,\, i\in \Z\}\subseteq \R_{\geq 0}$ is not dense.

\smallskip

For any perfectoid field $K$, we choose some element $\omega\in K^{\times}$ such that $|p|\leq |\omega|<1$. We define
$$K^{\flat\circ}:=\lim\limits_{\overleftarrow{x\mapsto \Phi(x)}}K^{\circ}/\omega.$$

Recall \cite[Lemma 3.2]{Scholze2012}.
\begin{lem}\label{lemaltdeskb}
\begin{points}
\item There exists a multiplicative homeomorphism $$\lim\limits_{\overleftarrow{x\mapsto x^p}} K^{\circ}\,\,\underrightarrow{\simeq} \,\,\lim\limits_{\overleftarrow{x\mapsto \Phi(x)}}K^{\circ}/\omega=K^{\flat\circ}$$ given by projection. Moreover, we have a map $$K^{\flat\circ}=\lim\limits_{\overleftarrow{x\mapsto x^p}} K^{\circ}=\{(x^{(0)},x^{(1)},\cdots)|\,\, x^{(i)}\in K^{\circ}, (x^{(i+1)})^p=x^i\}\to K^{\circ}$$ defined by $$x=(x^{(0)},x^{(1)},\cdots)\to x^{\#}:=x^{(0)}.$$ We may define a norm on $K^{\flat\circ}$ by $|x^{\#}|=|x|$ for all $x\in K^{\flat\circ}$.
\item 
The addition on $$K^{\flat\circ}=\{x:=(x^{(0)},x^{(1)},\cdots)|\,\, x^{(i)}\in K^{\circ}, (x^{(i+1)})^p=x^i\}$$ given by $(x+y)^i=\lim_{n\to \infty}(x^{(i+n)}+y^{(i+n)})^{p^n}.$
\item There exists an element $\omega^{\flat}\in \lim\limits_{\overleftarrow{x\mapsto x^p}} K^{\circ}$, satisfying $(\omega^{\flat})^{\#}=\omega$. Define
$$K^{\flat}:=K^{\flat\circ}[(\omega^{\flat})^{-1}].$$ Then norm $|\cdot|$ on $K^{\flat\circ}$ extends to a norm on $K^{\flat}$ which makes $K^{\flat\circ}$ to be the valuation ring of $K^{\flat}$. \item There exists a multiplicative homeomorphism $$K^{\flat}\underrightarrow{\simeq} \lim\limits_{\overleftarrow{x\mapsto x^p}}K.$$ Then $K^{\flat}$ is a perfectoid field of characteristic $p$. We have $|K^{b\times}|=|K^{\times}|$, $K^{\flat\circ}/\omega^{\flat}\simeq K^{\circ}/\omega$, and $K^{\flat\circ}/\mathfrak{m}^{\flat}\simeq K^{\circ}/\mathfrak{m},$ where $\mathfrak{m}$, \emph{resp.} $\mathfrak{m}^{\flat}$, is the maximal ideal of $K^{\circ}$, \emph{resp.} $K^{\flat\circ}$.
\item If $K$ is of characteristic $p$, then $K^{\flat}=K$.
\end{points}
\end{lem}
We note that (i) and (ii) of Lemma \ref{lemaltdeskb} implies that 
the definition of $K^{\flat\circ}$ is independent of $\omega$.

\medskip

We call $K^{\flat}$ the \emph{tilt} of $K.$
\begin{exe}The tilt of $\C_p$  is $\C_p^{\flat}=\widehat{\overline{\F_p((t))}}$.
\end{exe}

Then we have the following theorem, which was known by the classical work of Fontaine-Wintenberger \cite{Fontaine1979}
\begin{thm}\label{thmfwperfeq}
\begin{points}
\item Let $L$ be a finite extension of $K$. Then $L$ with its natural topology induced by $K$ is a perfectoid field.
\item The tilt functor $L\mapsto L^{\flat}$ induces an equivalence of categories between the category of finite extensions of $K$ and the category of finite extensions of $K^{\flat}$. This equivalence preserves degrees.
\end{points}
\end{thm}

\subsection{Almost mathematics}
Let $K$ be a perfectoid field and $\mathfrak{m}$ be the maximal ideal of $K^{\circ}$.

A $K^{\circ}$-module $M$ is said to be \emph{almost zero} if $\mathfrak{m}M=0.$
Define the category of almost $K^{\circ}$-modules as $$K^{\circ a}\text{-}\rm{mod} := K^{\circ}\text{-}\rm{mod} / ({\mathfrak{m}}\text{-}\text{torsion}).$$ We have a localization functor $M \mapsto M^a$ from $K^{\circ }$-mod to $K^{\circ a}$-mod, whose kernel is the thick subcategory of almost zero modules.

%

For two $K^{\circ a}$-modules $X,Y$, we define $\alHom(X,Y)=\Hom(X,Y)^a$.

\begin{pro}[\cite{Gabber2003}]\label{protensoral} The Category $K^{\circ a}$-mod is an abelian tensor category, where we define kernels, cokernels and tensor products in the unique way compatible with their definition in $K^{\circ}-\mod$, that is

$$ M^{a} \otimes N^a = (M \otimes N)^a$$

for any two $K^{\circ}$-modules $M,N$. For any $L,M,N \in K^{\circ a}$-mod there is a functorial isomorphism

$$ \Hom(L, \alHom(M,N))= \Hom(L\otimes M,N).$$
\end{pro}
This means that $K^{\circ a}$-mod  has all properties of the category of modules over a ring and thus one can define the notion of $K^{\circ a}$-algebra. Any $K^{\circ}$-algebra $R$ defines a $K^{\circ a}$-algebra $R^{a}$ as the tensor products are compatible. 
Moreover, localization also gives a functor from $R$-modules to $R^{a}$-modules.

\begin{pro}[\cite{Gabber2003}]\label{prorigadj}There exists a right adjoint functor
$$K^{\circ a}\text{-}\rm{mod} \to K^{\circ}\text{-}\rm{mod}: M\mapsto M_*:=\Hom_{K^{\circ a}}(K^{\circ a}, M).$$
to the localization functor $M\mapsto M^a$. The adjunction morphism $(M_*)^a\to M$ is an isomorphism. If $M$ is a $K^{\circ}$-module, then $(M^a)_*=\Hom(\mathfrak{m},M).$
\end{pro}

If $A$ is a $K^{\circ a}$-algebra, then $A_*$ has a natural structure as $K^{\circ}$-algebra and $(A^a) _* = A$. In particular, any $K^{\circ a}$-algebra comes via localization from a $K^{\circ}$-algebra.
Furthermore the functor $M\mapsto M_*$ induces a functor from $A$-modules to $A_*$-modules, and one can see also that all $A$-modules come via localization from $A_*$-modules. The category of $A$-modules is again an abelian tensor category, and all properties about the category of $K^{\circ a}$-modules stay true for the category of $A$-modules. 

Let $A$ be any $K^{\circ a}$-algebra.
As in \cite{Scholze2012}, a $A$-module $M$ is said to be flat if the functor $X \mapsto M \otimes _{A} X$ on $A$-modules is exact.

Denote by $\omega$ an element in $K^{\circ}$ satisfying $|p|\leq |\omega| <1$.
Let $A$ be a $K^a$-algebra, we say $A$ is $\omega$-adically complete if $A\simeq \lim\limits_{\longleftarrow}A/\omega^n.$

\subsection{Perfectoid algebras}
Fix a perfectoid field $K$ and an element $\omega\in K^{\circ}$ satisfying $|p|\leq |\omega| <1$.

\begin{defi}\begin{points}
\item A \emph{perfectoid $K$-algebra} is a Banach $K$-algebra $R$ such that the subset $R^{\circ}\subseteq R$ of powerbounded  elements is open and bounded, and the Frobenius morphism $\Phi: R^{\circ}/\omega \to R^{\circ}/\omega$ is surjective. Morphisms between perfectoid $K$-algebras are the continuous morphisms of $K$-algebras.
\item \emph{A perfectoid $K^{\circ a}$-algebra} is a $\omega$-adically complete flat $K^{\circ a}$-algebra $A$ on which Frobenius induces an isomorphism
$$\Phi: A/\omega^{1/p}\simeq A/\omega.$$ Morphism between perfectoid $K^{\circ a}$-algebras are the morphisms of $K^{\circ a}$-algebras.
\item\emph{ A perfectoid $K^{\circ a}/\omega$-algebra} is a flat $K^{\circ a}/\omega$-algebra $\bar{A}$ on which Frobenius induces an isomorphism
$$\Phi: \bar{A}/\omega^{1/p}\simeq \bar{A}.$$ Morphisms are the morphisms of $K^{\circ a}/\omega$-algebras.
\end{points}
\end{defi}

Let $K\text{-}\Perf$ denote the category of perfectoid $K$-algebras and similarly for $K^{\circ a}\text{-}\Perf$ and $K^{\circ a}/\omega\text{-}\Perf$. Let $K^{\flat}$ be the tilt of $K$ and $\omega^{\flat}$ is an element in $K^{\flat}$ satisfying $(\omega^{b})^{\#}=\omega$.

We recall \cite[Theorem 5.2]{Scholze2012}.
\begin{thm}\label{thmperfcateq}We have the following series of equivalences of categories:
$$K\text{-}\Perf\simeq K^{\circ a}\text{-}\Perf\simeq (K^{\circ a}/\omega)\text{-}\Perf=(K^{\flat a}/\omega^{b})\text{-}\Perf\simeq K^{\flat a}\text{-}\Perf\simeq K^{\flat}\text{-}\Perf.$$
\end{thm}
In other words, a perfectoid $K$-algebra, which is an object over the generic fibre, has a canonical extension to the almost integral level as a perfectoid $K^{\circ a}$-algebra, and perfectoid $K^{\circ a}$-algebras are determined by their reduction modulo $\omega$.

Let $R$ be a perfectoid $K^{\circ a}$-algebra, with $A=R^{\circ a}$. Define $$A^{\flat}:=\lim\limits_{\overleftarrow{\Phi}}A/\omega,$$
which we regard as a $K^{\flat \circ a}$-algebra via
$$K^{\flat \circ a}=(\lim\limits_{\overleftarrow{\Phi}}K^{\circ}/\omega)^a=\lim\limits_{\overleftarrow{\Phi}}(K^{\circ}/\omega)^a=\lim\limits_{\overleftarrow{\Phi}}K^{\circ a}/\omega,$$
and set $R^{\flat}=A^{b}_*[(\omega^{\flat})^{-1}].$

\begin{pro}\label{pro517}This defines a perfectoid $K^{\flat}$-algebra $R^{\flat}$ with corresponding perfectoid $K^{\flat \circ a}$-algebra $A^{\flat}$, and $R^{\flat}$ is the tilt of $R$. Moreover,
$$R^{\flat}=\lim_{\overleftarrow{x\mapsto x^p}}R, A^{\flat}_*=\lim\limits_{\overleftarrow{x\mapsto x^p}} A_*, \text{ and } A^{\flat}_*/\omega^{\flat}\simeq A_*/\omega. $$
In particular, we have a continuous multiplicative map $R^{\flat}=\lim\limits_{\overleftarrow{x\mapsto x^p}}R\to R$, $$x=(x^{(0)},x^{(1)},\cdots)\mapsto x^{\#}:=x^{(0)}.$$
\end{pro}

Then the equivalence $K\text{-}\Perf \to K^{\flat}\text{-}\Perf$ in Theorem \ref{thmperfcateq} is given by $R\mapsto R^{\flat}$.

\begin{pro}\label{prppolydisq}Let $R=K\langle T_1^{1/p^{\infty}},\cdots,T_n^{1/p^{\infty}}\rangle=K^{\circ}\widehat{[T_1^{1/p^{\infty}},\cdots,T_n^{1/p^{\infty}}]}[\omega^{-1}].$
Then $R$ is a perfectoid $K$-algebra, and its tilts $R^{\flat}$ is given by $K^{\flat}\langle T_1^{1/p^{\infty}},\cdots,T_n^{1/p^{\infty}}\rangle.$
\end{pro}

\subsection{Perfectoid spaces}
Fix a perfectoid field $K$ and an element $\omega\in K^{\circ}$ satisfying $|p|\leq |\omega| <1$. Let $K^{\flat}$ be the tilt of $K$ and $\omega^{\flat}$ is an element in $K^{\flat}$ satisfying $(\omega^{b})^{\#}=\omega$.

\begin{defi}An \emph{perfectoid  affinoid $K$-algebra} is an affinoid $K$-algebra $(R,R^{+})$, where $R$ is a perfectoid $K$-algebra, and $R^{+}\subseteq R^{\circ}$ is an open and integrally closed subring.
\end{defi}

\begin{pro}\label{protiltprealg}The association $(R,R^{+})\mapsto (R^{b},R^{\flat +})$, where $R^{\flat +}=\lim\limits_{\overleftarrow{x\to x^p}}R^{+}.$  defines an equivalence between the category of perfectoid affinoid $K$-algebras and the category of perfectoid affinoid $K^{\flat}$-algebras.
\end{pro}

%

\begin{thm}\label{thmtopiso}For any $x\in \Spa(R,R^+)$, one may define a point $x^{\flat}\in \Spa(R^{\flat},R^{\flat +})$ by setting $|f(x^{\flat})|:=|f^{\#}(x)|$ for $f\in R^{\flat}$.
This defines a homeomorphism $\rho:\Spa(R,R^+)\underrightarrow{\simeq} \Spa(R^{\flat},R^{\flat +})$ preserving rational subsets.
\end{thm}
Denote by $X:=\Spa(R,R^+)$ and $X^{\flat}:=\Spa(R^{\flat},R^{\flat +}).$
We note that in general the map $R^{\flat}\to R: f\to f^{\#}$ is not surjective. For any $f$ in $R$, $\rho_*f:=f\circ \rho^{-1}$ is a continuous function on $X^{\flat}$ but in general is not contained in $R^{\flat}.$

\smallskip

We have the following approximation lemma \cite[Corollary 6.7]{Scholze2012}.

\begin{lem}\label{lemappro}For any $f\in R$ and any $c\geq 0, \epsilon >0$, there exists $g_{c,\epsilon}\in R^{\flat}$ such that for all $x\in X$, we have
$$|f(x)-g_{c,\epsilon}^{\#}(x)|\leq |\omega|^{1-\epsilon}\max(|f(x)|,|\omega|^c).$$
\end{lem}

\begin{rem}Note that for $\epsilon<1$, the given estimate says in particular that for all $x\in X$, we have
$$\max\{|f(x)|, |\omega|^c\}=\max\{|g_{c,\epsilon}^{\#}(x)|,|\omega|^c\}.$$
\end{rem}

\begin{rem}\label{remappro}
Let $R:=K\langle x_1,\dots,x_N\rangle$ and $R^+:=R^{\circ}=K^{\circ}\langle x_1,\dots,x_N\rangle$. Then $R^{\flat}=K^{\flat}\langle x_1,\dots,x_N\rangle$ and $R^{\flat +}=K^{\flat\circ}\langle x_1,\dots,x_N\rangle$.

By Lemma \ref{lemappro}, for any $c \in \Z^+$, there exists an element $g_{c}\in K^{\flat \circ}\langle x_1^{1/p^{\infty}},\dots, x_N^{1/p^{\infty}}\rangle$ such that for all $x\in U_0^{\perf}$, we have
$$|H\circ \pi(x)-g_{c}^{\#}(x)|\leq |p|^{1/2}\max(|H\circ \pi(x)|,|p|^c).$$ There exists ${\ell}\in \N$ and an element $G_c\in K^{\flat \circ}[x_1^{1/p^{{\ell}}},\dots, x_N^{1/p^{{\ell}}}]$ such that $g_c-G_c\in t^{c+1} K^{\flat \circ}\langle x_1^{1/p^{\infty}},\dots, x_N^{1/p^{\infty}}\rangle$. It follows that for all $x\in U_0^{\perf}$, we have
$$|H\circ \pi(x)-G_{c}^{\#}(x)|\leq |p|^{1/2}\max(|H(x)|,|p|^c)=|p|^{1/2}\max(|G_{c}^{\#}(x)|,|p|^c),$$ and $G^{p^{\ell}}\in K^{\flat \circ}[x_1,\dots,x_N].$

Moreover, when $K^{\flat}=\widehat{\overline{F_p((t))}},$ we may make that $G_c\subseteq E^{\circ}[x_1^{1/p^{{\ell}}},\dots, x_N^{1/p^{{\ell}}}]$ where $E$ is a finite extension of $\overline{F_p}((t)).$

\end{rem}

\bigskip

Then we describe the structure sheaf $O_X$ on $X:=\Spa(R,R^+).$ Let $U=U(f_1,\cdots,f_n;g)\subseteq X$ be a rational subset. Equip $R[g^{-1}]$ with the topology making the image of $R^+[f_1/g,\cdots,f_n/g]\to R[g^{-1}]$ to be open and bounded. Let $R\langle f_1/g,\cdots,f_n/g \rangle$ be the completion of $R[g^{-1}]$ with respect to this topology. It equip with a subring $$R\langle f_1/g,\cdots,f_n/g\rangle^+\subseteq R\langle f_1/g,\cdots,f_n/g\rangle$$ which is the completion of integral closure of $R^+[f_1/g,\cdots,f_n/g]$.
By \cite[Proposition 1.3]{Huber1994}, the pair $(O_X(U),O_X^+(U)):=(R\langle f_1/g,\cdots,f_n/g\rangle, R\langle f_1/g,\cdots,f_n/g\rangle^+)$ depends only on the rational subset $U\subseteq X$ ( and not on the choice of $f_1,\cdots,f_n, g\in R$). The map
$$\Spa(O_X(U),O_X^+(U))\to \Spa(R,R^+)$$
is a homeomorphism onto $U$, preserving rational subsets. Moreover, $(O_X(U),O_X^+(U))$ is initial with respect to this property.

By \cite[Theorem 6.3]{Scholze2012}, we have the following
\begin{thm}\label{thmpersshesf}For any rational subset $U\subseteq X$, let $U^{\flat}:=\rho(U)\subseteq X^{\flat}.$
\begin{points}
\item The presheaves $O_X,O_{X^{\flat}}$ are sheaves.
\item For any rational subset $U\subseteq X$, the pair $(O_X(U),O_X^+(U))$ is a perfectoid affinoid $K$-algebra, which tilts to $(O_{X^{\flat}}(U^{\flat}),O_{X^{\flat}}^+(U^{\flat})).$
\end{points}
\end{thm}

The resulting spaces $\Spa(R,R^+)$ equipped with the two structure sheaves of topological rings $O_X,O_{X}^+$ are called \emph{affinoid perfectoid spaces} over $K$. The morphisms betweens the affinoid perfectoid spaces over $K$ are the morphisms induced by the morphisms between affinoid perfectoid $K$-algebras.

One defines \emph{perfectoid spaces} over $K$ to be the objects obtained by gluing affinoid perfectoid spaces. The morphisms betweens the perfectoid spaces are the morphisms glued by the morphisms between affinoid perfectoid spaces.

\medskip

We say that a perfectoid space $X^{\flat}$ over $K^{\flat}$ is the tilt of a perfectoid space $X$ over $K$ if there exists a functorial isomorphism $\Hom(\Spa(R^{\flat},R^{\flat +}),X^{\flat})=\Hom(\Spa(R,R^+),X)$ for all perfectoid affinoid $K$-algebras $(R,R^+)$ with tilts $(R^{\flat},R^{\flat +}).$

\begin{thm}\label{thmeqcatpersp}Any perfectoid space $X$ over $K$ admits a tilt $X^{\flat}$, unique up to isomorphism. This induces an equivalence between the category of perfectoid spaces over $K$ and the category of perfectoid spaces over $K^{\flat}$. The underly topological spaces of $X$ and $X^{\flat}$ are naturally identified by $\rho$. A perfectoid space $X$ is affinoid perfectoid if and only if its tilt $X^{\flat}$ is affinoid perfectoid. Finally, for any affinoid perfectoid subspace $U\subseteq X,$ the pair $(O_X(U),O_X^+(U))$ is a perfectoid affinoid $K$-algebra with tilt $(O_{X^{\flat}}(U^{\flat}),O^+_{X^{\flat}}(U^{\flat})).$
\end{thm}

For any morphism $F: X\to Y$ between perfectoid spaces over $K$, denote by $F^{\flat}:X^{\flat}\to Y^{\flat}$ the morphism between perfectoid spaces over $K^{\flat}$ induced by the equivalence of categories.

\subsection{Points in perfectoid spaces}
Fix a perfectoid field $K$ and an element $\omega\in K^{\circ}$ satisfying $|p|\leq |\omega| <1$.
Let $X$ be a perfectoid space over $K$.

For any point $x\in X$, let $K(x)$ be the residue field of $O_{X,x}$ and $K(x)^+\subseteq K(x)$ be the image of $O_{X,x}^+$.
By \cite[Proposition 2.25]{Scholze2012}, the $\omega$-adic completion of $O_{X,x}^+$ is equal to the $\omega$-adic completion $\widehat{K(x)}^+$ of $\widehat{K(x)}^+$. By \cite[Corollary 6.7]{Scholze2012}, $\widehat{K(x)}$ is a perfectoid field.

\begin{defi} An affinoid perfectoid field is a pair $(K,K^+)$ consisting of a perfectoid field and an open valuation subring $K^+\subseteq K.$
\end{defi}

Then $(\widehat{K(x)},\widehat{K(x)}^+)$ is an affinoid perfectoid field. Also note that affinoid perfectoid fields $(L,L^{+})$ for which $K\subseteq L$ are affinoid $K$-algebra.


\medskip

Then we have the following description of points \cite[Proposition 2.27]{Scholze2012}
\begin{pro}The points of $X$ are in bijection with maps $\iota: \Spa(L,L^+)\to X$ to affinoid fields $(L,L^+)$ such that the quotient field of the image of  $O_{X,x}$ in $L$ is dense, where $x$ is the image of $\Spa(L,L^+)$ in $X$.
\end{pro}

Any point $x\in X$ associates to a map $\iota: \Spa(\widehat{K(x)},\widehat{K(x)}^{+})\to X.$ By the equivalence of categories, the point $x^{\flat}\in X^{\flat}$ associates to a map $$\iota^{\flat}: \Spa(\widehat{K(x)}^{\flat},\widehat{K(x)}^{\flat +})\to X.$$ 

By \cite[Lemma 5.21]{Scholze2012}, $\Spa(\widehat{K(x)}^{\flat},\widehat{K(x)}^{\flat +})$ is an affinoid perfectoid field.
It follows that $\widehat{K^{\flat}(x^{\flat})}=(\widehat{K(x)})^{\flat}.$

In particular, we have the following
\begin{lem}\label{lembijripo}
For any point $x\in X$, $x$ is a $K$-point if and only if $x^{\flat}$ is a $K^{\flat}$-point in $X^{\flat}.$
\end{lem}

\section{Inverse limit of lifts of Frobenius}\label{sectioninvlim}
In this section, fix a perfectoid field $K$. Denote by $p>0$ the characteristic of the residue field $K^{\circ}/K^{\circ\circ}$ of $K$. 

Let $F:\P^N_{K}\rightarrow \P^N_{K}$ be a lift of Frobenius i.e. an endomorphism taking form $$F:[x_0:\cdots:x_N]\mapsto [x_0^q+p'P_0(x_0,\cdots,x_N):\cdots: x_N^q+p'P_N(x_0,\cdots,x_N)]$$ where $p'\in K^{\circ\circ}$, $q=p^s$ is a power of $p$, and $P_0,\cdots,P_N$ are homogeneous polynomials of degree $q$ in $K^{\circ}[x_0,\cdots,x_N].$
Let $\omega\in K^{\circ}$ be an element satisfying $\max\{|p'|,|p|\}\leq |\omega| <1$.

\subsection{Adic projective spaces}\label{subsectionadicproj}
At first, we define an adic space $\P^{N,\ad}_K$ which associates to the projective space $\P^{N}_K$.
In fact, by \cite[Theorem 2.22]{Scholze2012}, for any projective variety $X$ defined over $K$ with an integral model $\mathfrak{X}$ over $K^{\circ}$, we may associate an adic space $X^{\ad}$. But in this paper, we don't need the general theory and we define $\P^{N,\ad}_K$ in the following explicit way:

For any $i\in \{0,\cdots,N\}$, denote by
$$U_{i}^{\ad}:=\Spa(K\langle z_{i,0},\cdots,z_{i,i-1},z_{i,i+1},\cdots,z_{i,N}\rangle, K\langle z_{i,0},\cdots,z_{i,i-1},z_{i,i+1},\cdots,z_{i,N}\rangle^{\circ})$$ the unit balls.
Then we define $\P^{N,\ad}_K$ by gluing the unit balls $U_{i}^{\ad}$ together in the usual way:
For any $i\neq j$, $U_i^{\ad}\cap U_j^{\ad}=U(1,z_{i,0},\cdots,z_{i,i-1},z_{i,i+1},\cdots,z_{i,N}; z_{i,j})\subseteq U_i^{\ad}.$
On $U_i^{\ad}\cap U_j^{\ad}$, the transition map $\phi_{i,j}$ is defined by
$$\phi^*_{i,j}(z_{j,k})=z_{i,k}/z_{i,j} \text{ for } k\neq i,j; \text{ and } \phi^*_{i,j}(z_{j,i})=1/z_{i,j}.$$

Denote by  $R(\P^{N,\ad}_{K})$ the set of $K$-points in $\P^{N,\ad}_{K}$.

\begin{lem}\label{lemembrigid}There exists a natural embedding $\tau: \P^N_{K}(K)\hookrightarrow \P^{N,\ad}_{K}$. Moreover its image $\tau(\P^N_{K}(K))=R(\P^{N,\ad}_{K}).$
\end{lem}

\proof[Proof of Lemma \ref{lemembrigid}]
For any point $q\in \P^N_{K}(K)$,  there exists a finite extension $L$ of $K$, and a point $q'=[x_0:\dots:x_N]\in \P^N_{L}(L)$, such that $q$ is the image of that $q'$ under  the natural morphism $\pi^L_K:\P^N_L\to \P^N_K$ induced by the inclusion $K\hookrightarrow L.$ Indeed, $(\pi^L_K)^{-1}$ is exactly the Galois orbit of $q'.$ 
We may suppose that $\max\{|x_0|,\cdots,|x_N|\}= 1$ for all $j=0,\dots,N$. Denote by $I_q:=\{i|\,\, |x_i|=1\}.$ Observe that $I_q$ depends only on $q$. Pick $i\in I_q$, we define $\tau(q)\in U_i$ to be the point defined by  $f\to |f(x_1/x_i,\dots,x_n/x_i)|$, for all $f\in K\langle z_{i,0},\dots,z_{i,i-1},z_{i,i+1},\dots,z_{i,N}\rangle.$ Here $f(x_1/x_i,\dots,x_n/x_i)\in L$ depends on the choice of $q'$ in its Galois obit, but the value $|f(x_1/x_i,\dots,x_n/x_i)|$ depends only on $q$. Moreover we may check that the definition of $\tau(q)$ does not depend on the choice of $i\in I_q.$ Then $\tau$ is well defined. Moreover it is easy to check that $\tau$ is injective and $\tau(\P^N_{K}(K))\subseteq R(\P^{N,\ad}_{K})$. By \cite[6.1.2 Corollary 3]{Bosch1984}, the map $\tau: \P^N_{K}(K)\to R(\P^{N,\ad}_{K})$ is surjective.
\endproof

\subsection{Lifts of Frobenius on $\P^{N,\ad}_K$}\label{subsecliftfro} The endomorphism $F$ induces a natural endomorphism $F^{\ad}$ on $\P^{N,\ad}_K$. We define $F^{\ad}$ in the following explicit way.
For any $i=0,\dots,N$, $F^{\ad}|_{U_i^{\ad}}:U_i^{\ad}\to U_i^{\ad}$ is defined to be $$F^*(z_{i,j})=\frac{z_{i,j}^q+p'P_j(z_{i,0},\dots,z_{i,i-1},1,z_{i,i+1},\dots,z_{i,N})}{1+p'P_i(z_{i,0},\dots,z_{i,i-1},1,z_{i,i+1},\dots,z_{i,N})}$$ for all $j\neq i$.
We may write $$\frac{z_{i,j}^q+p'P_j(z_{i,0},\dots,z_{i,i-1},1,z_{i,i+1},\dots,z_{i,N})}{1+p'P_i(z_{i,0},\dots,z_{i,i-1},1,z_{i,i+1},\dots,z_{i,N})}=z_{i,j}^q+p'Q_{i,j}(z_{i,0},\dots,z_{i,i-1},z_{i,i+1},\dots,z_{i,N})$$
where $Q_{i,j}\in K^{\circ}\langle z_{i,0},\dots,z_{i,i-1},z_{i,i+1},\dots,z_{i,N}\rangle$.
For any $i\neq j$, we may check that $F_{i}^{\ad}(U_i^{\ad}\cap U_j^{\ad})\subseteq U_i^{\ad}\cap U_j^{\ad}$ and
$$F^{\ad}_i|_{U_i^{\ad}\cap U_j^{\ad}}=F^{\ad}_j|_{U_i^{\ad}\cap U_j^{\ad}}.$$

Then we may glue these $F_i^{\ad}$ to define $F^{\ad}: \P^{N,\ad}_K\to \P^{N,\ad}_K$.
Observe that we have the following commutative diagram:
\begin{equation*}
  \xymatrix{
    \P^N_K(K) \ar[d]_{F|_{\P^N_K(K)}} \ar[r]^\tau &\P^{N,\ad}_K\ar[d]^{F^{\ad}}\\
    \P^N_K(K) \ar[r]^{\tau^{\flat}} &\P^{N,\ad}_K}
\end{equation*}
Now we identify $\P^N_K(K)$ ( resp. $\P^N_{K^{\flat}}(K^{\flat}))$ to the image of $\tau$ ( resp. $\tau^{\flat}$) in $\P^{N,\ad}_K$ ( resp. $\P^{N,\ad}_{K^{\flat}}$).

\subsection{The inverse limit }
The inverse limit $\lim\limits_{\overleftarrow{F^{\ad}}}\P^{N,\ad}_K$ is the topological space $\{(x_0,x_1,\dots)\in (\P^{N,\ad}_K)^{\N}|\,\, F^{\ad}(x_i)=x_{i-1} \text{ for all } i\geq 1\}$ with the product topology. There exists a natural automorphism $T$ on $\lim\limits_{\overleftarrow{F^{\ad}}}\P^{N,\ad}_K$ defined by $$T: (x_0,x_1,\dots)\to (F^{\ad}(x_0),x_0,x_1,\dots).$$

The aim of this section is to construct a perfectoid space $(\P^{N}_K)^{\perf}$ with an automorphism $F^{\perf}$ such that the topological dynamical system $((\P^{N}_K)^{\perf}, F^{\perf})$ is isomorphic to $(\lim\limits_{\overleftarrow{F^{\ad}}}\P^{N,\ad}_K, T).$

Since $(F^{\ad})^{-1}(U_i^{\ad})\subseteq U_i^{\ad}$ for all $i=0,\dots,N,$ we have
$$\lim\limits_{\overleftarrow{F^{\ad}}}\P^{N,\ad}_K=\cup_{i=0}^N(\lim\limits_{\overleftarrow{F^{\ad}}}U_i^{\ad}).$$
Moreover we have $T(U_i^{\ad})\subseteq U_i^{\ad}.$
It follows that we only need to construct a perfectoid affinoid space $U_i^{\perf}$ with an automorphism $F_i^{\perf}$ such that the topological dynamical system $(U_i^{\perf}, F_i^{\perf})$ is isomorphic to $(\lim\limits_{\overleftarrow{F^{\ad}}}U_i^{\ad}, T|_{U_i^{\ad}})$ and check that they can be glued together.

Denote by $R_i^n:=K\langle z_{i,0}^{(n)},\dots,z_{i,i-1}^{(n)},z_{i,i+1}^{(n)},\dots,z_{i,N}^{(n)}\rangle$ for all $i=0,\dots,N$ and $n\geq 0$.
We identify $z^{(0)}_{i,j}$ and $z_{i,j}$ for $i\neq j.$ For every $n\geq 0$, we have an embedding $R_i^{n}\hookrightarrow R_i^{n+1}$ defined by
$$z_{i,j}^{(n)}\mapsto (z_{i,j}^{(n+1)})^q+p'Q_{i,j}(z_{i,0}^{(n+1)},\dots,z_{i,i-1}^{(n+1)},z_{i,i+1}^{(n+1)},\dots,z_{i,N}^{(n+1)})$$
where $Q_{i,j}$ is defined in Section \ref{subsecliftfro}.
Then we denote by
$$R_i:=K\langle z_{i,0}^{(\infty)},\dots,z_{i,i-1}^{(\infty)},z_{i,i+1}^{(\infty)},\dots,z_{i,N}^{(\infty)}\rangle$$
the completion of $\cup_{n=0}^{\infty}R_i^{n}.$ Denote by $\|\cdot\|$ the norm on $R_i$ induced by the norms on $R_i^n, n\geq 0$.

\begin{lem}\label{lemriperf}
For every $i=1,\dots, N$, $R_i$ is a perfectoid $K$-algebra with $$R_i^{\circ}=K^{\circ}\langle z_{i,0}^{(\infty)},\dots,z_{i,i-1}^{(\infty)},z_{i,i+1}^{(\infty)},\dots,z_{i,N}^{(\infty)}\rangle.$$ Its tilt is given by $R_i^{\flat}=K^{\flat}\langle z_{i,0}^{1/p^{\infty}},\dots,z_{i,i-1}^{1/p^{\infty}},z_{i,i+1}^{1/p^{\infty}},\dots,z_{i,N}^{1/p^{\infty}}\rangle$.
\end{lem}
\proof[Proof of Lemma \ref{lemriperf}]
Observe that $K^{\circ}\langle z_{i,0}^{(\infty)},\dots,z_{i,i-1}^{(\infty)},z_{i,i+1}^{(\infty)},\dots,z_{i,N}^{(\infty)}\rangle$ is the completion of $\cup_{n=0}^{\infty}(R_i^{n})^{\circ}.$ It is easy to check that $$K^{\circ}\langle z_{i,0}^{(\infty)},\dots,z_{i,i-1}^{(\infty)},z_{i,i+1}^{(\infty)},\dots,z_{i,N}^{(\infty)}\rangle\subseteq R_i^{\circ}.$$

For any $f\in R_i$, there exists a sequence $f_n\in R_i^{n}$ such that $f_n\to f$ as $n\to \infty$.
There exists $M\geq 0$, such that for all $m,n\geq M,$ $\|f_n-f_m\|\leq 1$.
It follows that $f_n-f_M\in (R_{i}^n)^{\circ}$ for all $n\geq M$. Then $f-f_M\in \widehat{\cup_{n=0}^{\infty}(R_i^{n})^{\circ}}.$
If $\|f_M\|\leq 1$, we have $f_M\in (R_i^M)^{\circ}$ and then $f\in \widehat{\cup_{n=0}^{\infty}(R_i^{n})^{\circ}}.$
If $\|f_M\|> 1$, we have $\|f^n\|=\|f_M^n\|\to \infty$ as $n\to \infty$. Then $f$ is not power bounded. It follows that
$$R_i^{\circ}\subseteq K^{\circ}\langle z_{i,0}^{(\infty)},\dots,z_{i,i-1}^{(\infty)},z_{i,i+1}^{(\infty)},\dots,z_{i,N}^{(\infty)}\rangle.$$

It follows that $R_i^{\circ}$ is open and bounded.

We have $R_i^{\circ}/\omega=(K^{\circ}/\omega)\langle z_{i,0}^{(\infty)},\dots,z_{i,i-1}^{(\infty)},z_{i,i+1}^{(\infty)},\dots,z_{i,N}^{(\infty)}\rangle$ is the completion of $\cup_{n=0}^{\infty}R_i^n/\omega.$ The embedding $R_i^n/\omega\to R_i^n/\omega$ is given by
$$z_{i,j}^{(n)} \mod \omega\mapsto (z_{i,j}^{(n+1)})^q+p'Q_{i,j}(z_{i,0}^{(n+1)},\dots,z_{i,i-1}^{(n+1)},z_{i,i+1}^{(n+1)},\dots,z_{i,N}^{(n+1)})\mod \omega $$$$=(z_{i,j}^{(n+1)})^q \mod \omega.$$
It follows that $$R_i^{\circ}/\omega=(K^{\circ}/\omega)\langle z_{i,0}^{1/p^{\infty}},\dots,z_{i,i-1}^{1/p^{\infty}},z_{i,i+1}^{1/p^{\infty}},\dots,z_{i,N}^{1/p^{\infty}}\rangle.$$
Then the Frobenius morphism $\Phi:R_i^{\circ}/\omega\to R_i^{\circ}/\omega$ is surjective. It follows that $R_i$ is a perfectoid $K$-algebra.

By Proposition \ref{prppolydisq} and the categorical equivalence in Theorem \ref{thmperfcateq}, we have $R_i^{\flat}=K^{\flat}\langle z_{i,0}^{1/p^{\infty}},\dots,z_{i,i-1}^{1/p^{\infty}},z_{i,i+1}^{1/p^{\infty}},\dots,z_{i,N}^{1/p^{\infty}}\rangle.$
\endproof

We define $U_i^{\perf}:=\Spa(R_i,R_i^{\circ})$ and $F^{\perf}_i:U_i^{\perf}\to U_i^{\perf}$ the map induced by the morphism $R_i\to R_i$ defined by  $$z_{i,j}^{(n)}\to z_{i,j}^{(n-1)} \text{ for all } n\geq 1 \text{ and } z_{i,j}^{(0)}\to (z^{(0)}_{i,j})^q+p'Q_{i,j}(z_{i,0}^{(0)},\dots,z_{i,i-1}^{(0)},z_{i,i+1}^{(0)},\dots,z_{i,N}^{(0)}).$$

Then we define $(\P^{N}_K)^{\perf}$ by gluing $U_{i}^{\perf}$ together in the usual way:
For any $i\neq j$, $U_i^{\perf}\cap U_j^{\perf}=U(1,z_{i,0}^{(0)},\dots,z_{i,i-1}^{(0)},z_{i,i+1}^{(0)},\dots,z_{i,N}^{(0)}; z_{i,j}^{(0)})\subseteq U_i^{\perf}.$
On $U_i^{\perf}\cap U_j^{\perf}$, the transition map $\phi_{i,j}$ is defined to be
$$(\phi^{\perf}_{i,j})^*(z_{j,k}^{(n)})=z_{i,k}^{(n)}/z_{i,j}^{(n)} \text{ for } k\neq i,j; \text{ and } (\phi^{\perf}_{i,j})^*(z_{j,i}^{(n)})=1/z_{i,j}^{(n)}.$$
It is easy to check that for all $i\neq j$, $$F^{\perf}_i(U_i^{\perf}\cap U_j^{\perf})\subseteq U_i^{\perf}\cap U_j^{\perf}$$ and $F_i^{\perf}=F_j^{\perf}$ on $U_i^{\perf}\cap U_j^{\perf}.$ Then we define $F^{\perf}$ by gluing $F_i^{\perf}$ for $i=0,\dots,N.$

Then we have the following
\begin{thm}\label{thmdysyisoinvperf}There exists a natural homeomorphism $\psi:(\P^{N}_K)^{\perf}\to \lim\limits_{\overleftarrow{F^{\ad}}}\P^{N,\ad}_K$ makes the following diagram to be commutative:
\begin{equation*}
  \xymatrix{
    \P^{N,\perf}_K \ar[d]_{F^{\perf}} \ar[r]^\psi &\lim\limits_{\overleftarrow{F^{\ad}}}\P^{N,\ad}_K\ar[d]^{T}\\
    \P^{N,\perf}_K \ar[r]^{\psi} &\lim\limits_{\overleftarrow{F^{\ad}}}\P^{N,\ad}_K}
\end{equation*}
In other words, the topological dynamical systems $(\P^{N,\perf}_K, F^{\perf})$ and $(\lim\limits_{\overleftarrow{F^{\ad}}}\P^{N,\ad}_K,T)$ are isomorphic by $\psi.$

Moreover a point $x\in \P^{N,\perf}_K$ whose image $\psi(x)=(x_0,x_1,\dots)$ is a $K$-point if and only if for every $n\geq 0$, $x_n$ is a $K$-point.
\end{thm}

\proof[Proof of Theorem \ref{thmdysyisoinvperf}]Denote by $B_i:=\cup_{n=0}^{\infty}R_i^n$. We have $B_i^{\circ}=\cup_{n=0}^{\infty}R_i^{n\circ}.$ Then $\Spa(B,B^{\circ})$ is an affinoid space and we have $R_i=\widehat{B_i}$, $R_i^{\circ}=\widehat{B_i^{\circ}}.$ By Proposition \ref{procompletionaffi}, then natural morphism $\mu_i:\Spa(R_i,R_i^{\circ})\to \Spa(B_i,B_i^{\circ})$ is a homeomorphism.

\smallskip

Denote by $\psi_i^n: U^{\perf}_i\to U^{\ad}_i$ the map induced by the morphism $$K\langle z_{i,0},\dots,z_{i,i-1},z_{i,i+1},\dots,z_{i,N}\rangle\to R^{n}_i\subseteq B_i\subseteq R_i$$ by sending $z_{i,j}\to z_{i,j}^{(n)}$. It is easy to check that $\psi_i^n$ could be glued to a map $\psi^n:\P^{N,\perf}_K\to \P^{N,\ad}_K.$

Since $F^{\ad}\circ\psi^{n+1}=\psi^{n}$ for all $n\geq 0$, it induces a map $$\psi:=\lim\limits_{\overleftarrow{n}}\psi^n: \P^{N,\perf}_K\to \lim\limits_{\overleftarrow{F^{\ad}}}\P^{N,\ad}_K.$$
By checking in the affinoid spaces $U_i^{\perf}$'s, it is easy to check that $T\circ\psi=\psi\circ F^{\perf}.$

So we only need to show that $\psi$ is a homeomorphism. We only need to show it in $U_i.$ Denote by $$\psi_i:=\psi|_{U^{\perf}_i}=\lim\limits_{\overleftarrow{n}}\psi^n_i.$$

Now we define a morphism $\theta_i: \lim\limits_{\overleftarrow{F^{\ad}}}U_i^{\ad}\to \Spa(B_i,B_i^{\circ})$ as the following.
Let $(x_0,x_1,\dots)$ be a point in $\lim\limits_{\overleftarrow{F^{\ad}}}U_i^{\ad}$. For every $n\geq 0$, we identify $U_i\to \Spa(R_i^n,R_i^{n\circ})$ by $z_{i,j}\to z_{i,j}^{(n)}.$
Then $x_n$ defines a valuation on $R_i^n$ with valuation group
$\Gamma_n:=\{|f(x_n)||\,\,f\in R_i^n\}$. Moreover, for any ${\ell}\geq n$, and $f\in R_i^n$, we have $|f(x_l)|=|f(x_n)|.$
Then we define $\theta_i((x_0,x_1,\dots))$ to be the natural valuation $B_i=\cup_{n=0}^{\infty}R_i^n\to \cup_{n=0}^{\infty}\Gamma_n$ by gluing all the valuations $x_n$ on $R^n_i.$ Since all the rational subset of $\Spa(B_i,B_i^{\circ})$ are defined over some $R_i^n$, it is easy to check that $\theta$ is continuous. It is easy to check that $\psi_i\circ(\mu^{-1}_i\circ\theta_i)=\id$ and $(\mu_i^{-1}\circ\theta_i)\circ\psi_i=\id$. It follows that $\psi_i$ is a homeomorphism.

Let $x$ be a $K$-point in $U_i^{\perf}$ and $\psi_i(x)=(x_0,x_1,\dots).$ 
For any $n\geq 0$,
we have $$K\subseteq R_i^n/\{|f(x_n)|=0|\,\, f\in R^n_i\}\subseteq R_i/\{|f(x)|=0|\,\, f\in R_i\}=K.$$
It follows that $x_n$ is a $K$-point.

Let $x$ be a point in $U_i^{\perf}$ and set $\psi_i(x)=(x_0,x_1,\dots).$ We suppose that all $x_n$ are $K$-points. Then $m_i^n:=\{|f(x_n)|=0|\,\, f\in R^n_i\}$ is a maximal ideal in $R^n_i$ and 
$R^n_i/m_i^n=K.$ 
 The valuation $R^n_i/m_i^n\to \R$ induced by $x_i$ is  the norm on $K$.

There exists a continuous morphism $\cup_{n=0}^{\infty}R^{n}_i\to K$ by gluing the morphisms $R_i^n\to R_i^n/m_i^n\hookrightarrow K$. We can extend this morphism to a continuous morphism $g:R_i^n=\widehat{\cup_{n=0}^{\infty}R^{n}_i}\to K$. The valuation $f\to |g(f)|$  defines a point $y\in U_i^{\perf}$. Then $y$ is a $K$-point.
Observe that for all $f\in R_i^n$, $|f(y)|=|f(x_i)|$. Then we have $\psi(y)=(x_0,x_1,\dots)=\psi(x).$ Then $y=x$ and so,  $x$ is a $K$-point.
\endproof

For every $i=0,\dots,N$, the embedding
$K\langle z_{i,0},\dots,z_{i,i-1},z_{i,i+1},\dots,z_{i,N}\rangle\subseteq R_i$
induces a map $U_i^{\perf}\to U_i^{\ad}.$ We define $\pi: \P_{K}^{N,\perf}\to \P_{K}^{N,\ad}$ by gluing these maps.
It is easy to check that $$F^{\ad}\circ\pi^{\flat}=\pi^{\flat}\circ F^{\perf}.$$
For any point $x\in \P^{N,\perf}_{K}$ with $\psi(x)=(x_0,x_1,\dots)$,  we have $\pi(x)=x_0$.

\subsection{Passing to the tilt}
Denote by $U_i^{\flat,\perf}:=\Spa(R_i^{\flat},R_i^{b\circ})$ and $\Phi_i^{s,\perf}:U_i^{\flat,\perf}\to U_i^{\flat,\perf}$ the $s$-th power of the Frobenius i.e. the map induced by the morphism $$R_i^{\flat}\to R_i^{\flat} : f\to f^q.$$

We define $(\P^{N}_{K^{\flat}})^{\perf}$ by gluing $U_{i}^{\flat,\perf}$ together in the usual way:
For any $i\neq j$, $U_i^{\flat,\perf}\cap U_j^{\flat,\perf}=U(1,z_{i,0},\dots,z_{i,i-1},z_{i,i+1},\dots,z_{i,N}; z_{i,j})\subseteq U_i^{\flat,\perf}.$

On $U_i^{\flat,\perf}\cap U_j^{\flat,\perf}$, the transition map $\phi_{i,j}^{\flat}$ is defined to be
$$(\phi^{\flat,\perf}_{i,j})^*(z_{j,k}^{1/p^n})=z_{i,k}^{1/p^n}/z_{i,j}^{1/p^n} \text{ for } k\neq i,j; \text{ and } (\phi^{\flat,\perf}_{i,j})^*(z_{j,i}^{1/p^{n}})=1/z_{i,j}^{1/p^{n}}.$$
By reducing modulo $\omega$ and the categorical equivalence in Theorem \ref{thmperfcateq}, we see that $\phi_{i,j}^{\flat,\perf}=(\phi_{i,j}^{\perf})^{\flat}$ and $\Phi_i^{s,\perf}=(F^{\perf}_i)^{\flat}.$ It follows that $(\P^{N}_{K^{\flat}})^{\perf}=((\P^{N}_{K})^{\perf})^{\flat}$ and we can define $\Phi^{s,\perf}$ by gluing $\Phi_i^{s,\perf}$ together. Moreover, we have $\Phi^{s,\perf}=(F^{\perf})^{\flat}.$
Then we have the following
\begin{thm}\label{thmperfisdy}The following diagram is commutative:
\begin{equation*}
  \xymatrix{
    \P^{N,\perf}_K \ar[d]_{F^{\perf}} \ar[r]^\rho &\P^{N,\perf}_{K^{\flat}}\ar[d]^{\Phi^{s,\perf}}\\
    \P^{N,\perf}_K \ar[r]^{\rho} &\P^{N,\perf}_{K^{\flat}}}
\end{equation*}
In other words, the topological dynamical systems $(\P^{N,\perf}_K, F^{\perf})$ and $(\P^{N,\perf}_{K^{\flat}},\Phi^s)$ are isomorphic by $\rho.$
\end{thm}

For any $i\in \{0,\dots,N\}$, denote by
$$U_{i}^{\flat,\ad}:=\Spa(K^{\flat}\langle z_{i,0},\dots,z_{i,i-1},z_{i,i+1},\dots,z_{i,N}\rangle, K^{\flat}\langle z_{i,0},\dots,z_{i,i-1},z_{i,i+1},\dots,z_{i,N}\rangle^{\circ}).$$
As in Section \ref{subsectionadicproj}, we define $\P^{N,\ad}_{K^{\flat}}$ by gluing $U_i^{\flat,\ad}$, $i=0,\dots,N$. Denote by $\phi_i^{s,\ad}$ the $s$-th power of the Frobenius on $U_i$ i.e. the map induced by the morphism $f\to f^q$ on $K^{\flat}\langle z_{i,0},\dots,z_{i,i-1},z_{i,i+1},\dots,z_{i,N}\rangle.$

By Lemma \ref{lemembrigid}, we have a natural embedding $\tau^{\flat}: \P^N_{K^{\flat}}(K^{\flat})\hookrightarrow \P^{N,\ad}_{K^{\flat}}$.
Then $\tau(\P^N_{K}(K))=R(\P^{N,\ad}_{K})$ and we have the following commutative diagram
\begin{equation*}
  \xymatrix{
    \P^N_{K^{\flat}}(K^{\flat}) \ar[d]_{\Phi^s|_{\P^N_K(K)}} \ar[r]^\tau &\P^{N,\ad}_{K^{\flat}}\ar[d]^{\Phi^{s,\ad}}\\
    \P^N_{K^{\flat}}(K^{\flat}) \ar[r]^\tau &\P^{N,\ad}_{K^{\flat}}}
\end{equation*}
where $\Phi^s$ is the $s$-th power of the Frobenius on $P^N_{K^{\flat}}$.

For every $i=0,\dots,N$, the embedding
$K^{\flat}\langle z_{i,0},\dots,z_{i,i-1},z_{i,i+1},\dots,z_{i,N}\rangle\subseteq R_i^{\flat}$
induces a map $U_i^{\flat,\perf}\to U_i^{\flat,\ad}.$ We define $\pi^{\flat}: \P_{K^{\flat}}^{N,\perf}\to \P_{K^{\flat}}^{N,\ad}$ by gluing these maps.
It is easy to check that $$\Phi^{s,\ad}\circ\pi^{\flat}=\pi^{\flat}\circ\Phi^{s,\perf}   \eqno(1).$$
By \cite[Theorem 8.5]{Scholze2012}, $\pi^{\flat}$ is a homeomorphism.
Moreover, we have the following
\begin{lem}\label{lempibrigid}
The map $\pi^{\flat}$  induces a bijection between $R(\P^{N,\perf}_{K^{\flat}})$ and $R(\P^{N,\ad}_{K^{\flat}}).$
\end{lem}

\proof[Proof of Lemma \ref{lempibrigid}]It is clear that if $x$ is a $K^{\flat}$-point then $\pi^{\flat}(x)$ is a $K^{\flat}$-point.

Now we suppose that $\pi^{\flat}(x)$ is a $K^{\flat}$-point. We suppose that $x$ is contained in $U_i^{\flat,\perf}$ and then $x_0:=\pi^{\flat}(x)\in U_i^{\flat,\ad}.$
Since $x_0$ is a $K^{\flat}$-point, it defines a morphism $$g_0:R^{\flat,0}_{i}:=K^{\flat}\langle z_{i,0},\dots,z_{i,i-1},z_{i,i+1},\dots,z_{i,N}\rangle\to K^{\flat}.$$  It follows that the Frobenius map $f\to f^p$ on $K^{\flat}$ is a field automorphism. For any $f\in R^{\flat,n}_i:=K^{\flat}\langle z_{i,0}^{1/p^n},\dots,z_{i,i-1}^{1/p^n},z_{i,i+1}^{1/p^n},\dots,z_{i,N}^{1/p^n}\rangle$, we have $f^{p^n}\in R^{\flat,0}_{i}$. Then the morphism $g_0$ extends to a morphism  $g_n:R^{\flat,n}_i\to K^{\flat}$ by sending $f$ to $(g_0(f^{p^n}))^{1/p^n}.$ We glue $g_n$ to define a continuous morphism $\cup_{n=0}^{\infty}R^{\flat,n}_i\to K^{\flat}$ and then extend it to a continuous morphism $g: R^{\flat}_i=\widehat{(\cup_{n=0}^{\infty}R^{\flat,n}_i)}\to K^{\flat}$. Then $g$ induces a $K^{\flat}$-point $y\in U_i^{\flat,\perf}.$ Since $\pi^{\flat}(y)=x_0=\pi^{\flat}(x)$, we have $y=x$. Then $x$ is a $K^{\flat}$-point.
\endproof

\section{Periodic points}\label{sectionperiodic}
In this section, we denote by $K=\C_p$. Then $K$ is a perfectoid field and $K^{\flat}$ is the completion of the algebraical closure of $\F_p((t)).$ We may suppose that $|p|=|t|=p^{-1}.$

Let $F:\P^N_{K}\rightarrow \P^N_{K}$ be an endomorphism taking form $$F:[x_0:\dots:x_N]\mapsto [x_0^q+p'P_0(x_0,\dots,x_N):\dots: x_N^q+p'P_N(x_0,\dots,x_N)]$$ where $p'\in K^{\circ\circ}$, $q$ is a power of $p$, and $P_0,\dots,P_N$ are homogeneous polynomials of degree $q$ in $K^{\circ}[x_0,\dots,x_N].$
The aim of this section is to study the periodic points of $F$. In particular, we prove Theorem \ref{thmdvtperlf} and Theorem \ref{thmmmdlfper}.

\medskip

Recall that $\Per$ is the set of periodic closed points in $\P^N_{K}$.

Let $V$ be any irreducible subvariety of $\P^N_{K}$.
Suppose that $V$ is defined by the equations $H_j(x_0,\dots,x_N)=0$, $j=1,\dots,m$ where $H_j$ are homogenous polynomials. We may suppose that $\|H_j\|=1$ for all $j=1,\dots,m$. For any $i=0,\dots,N$, denote by
$$V_i^{\ad}:=\{x\in U_i^{\ad}|\,\, |H_{i,j}(x)|=0 , j=1,\dots,m\}$$ where $H_{i,j}:=H(z_{i,0},\dots,z_{i,i-1},1,z_{i,i+1},\dots,z_{i,N})$.
Observe that $\|H_{i,j}\|=1$.

 Set $R(V^{\ad}_i):=R(\P^{N,\ad}_K)\cap V_i^{\ad}$, $V^{\ad}:=\cup_{i=0}^NV_i^{\ad}$ and $R(V^{\ad}):=R(\P^{N,\ad}_K)\cap V^{\ad}$.
We have $\tau(R(V))=R(V^{\ad}).$

Observe that for all points $x\in R(U_i^{\ad})$, we have $d(x,V)=\max\{|H_{i,j}(x)|\}.$

\subsection{Passing to the reduction}
Since $K$ is algebraically closed, we have $\P^N_K(K)=\P^N_K(K).$ Denote by $k=\overline{\F_p}$, we have $k=K^{\circ}/K^{\circ\circ}.$
At first, there exists a reduction map $${\red}:\P^N_K(K)\to \P^N_k(k)$$ defined by the following:
For any point $x\in \P^N_K(K)$, we may write it as $x=[x_0:\dots:x_N]$ where $x_i\in K^{\circ}, i=0,\dots,N$ and $\max\{|x_i|,i=0,\dots,N\}=1.$ Then we define ${\red}(x)=[\overline{x_0}:\dots:\overline{x_N}]$ where $\overline{x_i}$ is the image of $x_i$ in $k=K^{\circ}/K^{\circ\circ}.$ Observe that ${\red}\circ F=\overline{\Phi}^s\circ \red$, where $\overline{\phi}$ is the Frobenius on $\P^N_k$.
For every point $y\in \P^N_k(k)$, there exists $m>0$ such that $\overline{\Phi}^{sm}(y)=y.$  Then we have $D_y:=\red^{-1}(y)\simeq (K^{\circ\circ})^N$ is a polydisc fixed by $F$.
Since $F^m|_{D_y}$ is attracting, $D_y\cap \Per$ has exactly one point. It follows that $\red$ induces a bijection between $\Per$ and $\P^N_k(k).$

 Similarly, we can define the reduction map $\red^{\flat}:\P^N_{K^{\flat}}(K^{\flat})\to \P^N_k(k)$. This map induces a bijection between $\Per^{\flat}$ and $\P^N_k(k)$ where $\Per^{\flat}$ is the set of $\Phi^s$-periodic closed points of $\P^N_{K^{\flat}}$.

Since $k$ is a subfield of $K^{\flat}$, there exists an embedding $\eta: \P^N_k(k)\hookrightarrow \P^{N}_{K^{\flat}}(K^{\flat}).$ Observe that the image $\eta(\P^N_k(k))$ is exactly $\Per^{\flat}$. Moreover we have $\red^{\flat}\circ\eta=\id.$
We may check that the map $$\phi:=\eta\circ\red:\Per\to \Per^{\flat}$$ is a bijection satisfying $\Phi^s\circ\phi=\phi\circ F$.

\subsection{Passing to the tilt}
Denote by $\Per^{\ad}=\tau(\Per).$ It is exactly the set of periodic $K$-points in $\P_{K}^{N,\ad}.$ For any point $x\in \Per^{\ad}$, denote by $n>0$ a period of $x$ under $F^{\ad}.$ We define a map $\chi:\Per^{\ad}\to \lim\limits_{\overleftarrow{F^{\ad}}}\P^{N,\ad}_{K}$ by sending $x$ to $(x_0,x_1,\dots)$ where
$x_i=(F^{\an})^{kn-i}(x)$ where $kn\geq i$. We note that $\chi(x)$ does not depend on the choice of $n$ and $k$.
Since $\pi\circ\chi=\id$,
$\chi$ is injective. We have that $\chi(\Per^{\ad})$ is exactly the set of $\Per_T$, where $\Per_T$ is the set of points $(x_0,x_1,\dots)\in \lim\limits_{\overleftarrow{F^{\ad}}}\P^{N,\ad}_{K}$ which is periodic under $T$ such that every $x_n$ is a $K$-point.

Denote by $\Per^{\flat,\ad}$ the set of $K^{\flat}$-points in $\P^{N,\ad}_{K^{\flat}}$ which are periodic under $\Phi^{s,\ad}.$
By applying Lemma \ref{lemembrigid} over $K^{\flat}$, there exists a bijection $\tau^{\flat}: {\P^N_{K^{\flat}}(K^{\flat})}\to R(\P_{K^{\flat}}^{N,\ad})$ and we have
$$\tau^{\flat}\circ\Phi^s=\Phi^{s,\ad}\circ\tau^{\flat}.$$ It follows that $\tau^{\flat}$ induces a bijection between $\Per^{\flat,\ad}$ and the set $\Per^{\flat}$ of $\Phi^s$-periodic points in  $R(\P^{N}_{K^{\flat}})=\P^{N}_{K^{\flat}}(K^{\flat})$.

By Theorem \ref{thmdysyisoinvperf}, Theorem \ref{thmperfisdy}, Equation (1) and Lemma \ref{lempibrigid}, the map
$$\iota:=\pi^{\flat}\circ\rho\circ\psi^{-1}\circ\chi:\Per^{\ad}\to \Per^{\flat,\ad}$$ is bijective.

Denote by $\Per^{\ad}_i:=\Per^{\ad}\cap U_i^{\ad}$ and $\Per^{\flat,\ad}_i:=\Per^{\flat,\ad}\cap U_i^{\flat,\ad}$ for every $i=0,\dots,N$.
Then we have $\iota(\Per^{\ad}_i)=\Per^{\flat,\ad}_i.$

Observe that for every point $x\in \Per$, we have $\red(x)=\red^{\flat}\circ\iota\circ\tau(x).$ Then we have $$\phi=\eta\circ \red=\eta\circ\red^{\flat}\circ\iota\circ\tau=\iota\circ\tau$$ on $\Per.$

\subsection{Proof of Theorem \ref{thmdvtperlf}}
We only need to show this theorem for the periodic points in $U_i^{\ad}$ for all $i=0,\dots, N.$ Without the loss of generality, we only need to show that there exists $\delta>0$ such that for all $x\in \Per\cap U_0^{\ad}$, either $d(x,V)>\delta$ or $x\in V.$

At first, we prove our theorem for hypersurfaces.
\begin{lem}\label{lemvtperhyper}Let $H\in K[x_1,\dots,x_N]$ be a polynomial. Then there exists $\varepsilon>0$, such that for all $x\in \Per\cap U_0^{\ad}$, either $|H(x)|>\varepsilon$ or $H(x)=0.$
\end{lem}

By this lemma, for any $H_{0,j}$, $j=1,\dots, m$, we have $\varepsilon_j>0$ such that for all $x\in \Per\cap U_0^{\ad}$, either $|H_{0,j}(x)|>\varepsilon_j$ or $H_{0,j}(x)=0.$ Set $\delta:=\min_{1\leq j\leq m}\{\varepsilon_j\}$. Let $x$ be a point in $\Per\cap U_0^{\ad}$ satisfying $d(x,V)\leq \delta$. Then for all $j=1,\dots, m$, we have $H_{0,j}(x)=0.$ It follows that $x\in V.$

\medskip

We only need to proof Lemma \ref{lemvtperhyper}.

To do this, we need the following Lemma

\begin{lem}\label{lemtotallyram}Let $E/\overline{F_p}((t))$ be a finite extension. Then $E=\overline{F_p}((u))$ for some $u\in E$ satisfying $|u|=|t|^{1/[E:\overline{F_p}((t))]}$.
\end{lem}
\proof[Proof of Lemma \ref{lemtotallyram}]
Observe that $E$ is a discrete valuation field.

Since $\overline{F_p}$ is algebraically closed, the extension $E/\overline{F_p}((t))$ is totally ramified. It follows that $E=\overline{F_p}((t))(u)$ where the minimal polynomial of $u$ over $\overline{F_p}((t))$ is an Eisenstein polynomial. It follows that $|u|=|t|^{1/[E:\overline{F_p}((t))]}$ and $uE^{\circ}$ is the maximal ideal of $E^{\circ}$. For every $f\in E^{\circ}$, $f$ can be written as $\sum_{i\geq 0}a_iu^i$ where $a_i\in \overline{F_p}$ for all $i\geq 0$. Then we conclude our proof.
\endproof

\begin{lem}\label{lemgoodrepapro}For any polynomial $G\in K^{\flat\circ}[x_1,\dots,x_N]$ and $\varepsilon>0$, there exists a polynomial $G_{\varepsilon}\in K^{\flat\circ}[x_1,\dots,x_N]$ satisfying $\deg G_{\varepsilon}\leq \deg G$, $\|G-G_{\varepsilon}\|<\varepsilon$ (resp. $\leq \varepsilon$)and $G_{\varepsilon}$
takes the form $G_{\varepsilon}=\sum_{i\geq 0}^{m}u^ig_i$ where $g_i\in \overline{\F_p}[x_1,\dots,x_N]$, $u\in \overline{\F_p((t))}^{\circ}$ with norm $|u|=|t|^{1/[\overline{\F_p}((t))(u):\overline{\F_p}((t))]}$ and $|u|^m\geq \varepsilon$ (resp. $>\varepsilon$).
\end{lem}
\proof[Proof of Lemma \ref{lemgoodrepapro}]By Lemma \ref{lemtotallyram}, there exists $u\in \overline{F_p((t))}^{\circ}$ with norm $|u|=|t|^{1/[\overline{\F_p}((t))(u):\overline{\F_p}((t))]}$ and $H\in \overline{F_p}((t))(u)[x_1,\dots,x_N]$ such that $\deg H\leq \deg G$, $\|G-H\|<\varepsilon$ and $H$
takes form $H=\sum_{i\geq 0}^{\infty}u^ig_i$ where $g_i\in \overline{\F_p}[x_1,\dots,x_N]$. Let $m$ be the largest integer such that $|u|^m\geq \varepsilon$ (resp. $>\varepsilon$). Set $G_{\varepsilon}=\sum_{i\geq 0}^{m}u^ig_i$ then we conclude our proof.
\endproof
\proof[Proof of Lemma \ref{lemvtperhyper}]We may suppose that $H\neq 0$ and $\|H\|=1.$

By Remark \ref{remappro}, for any $c \in \Z^+$, there exists ${\ell}\in \N$ and an element $G_c\in K^{\flat \circ}[x_1^{1/p^{{\ell}}},\dots, x_N^{1/p^{{\ell}}}]$ such that for all $x\in U_0^{\perf}$, we have
$$|H\circ \pi(x)-G_{c}^{\#}(x)|\leq |p|^{1/2}\max(|H(x)|,|p|^c)=|p|^{1/2}\max(|G_{c}^{\#}(x)|,|p|^c),$$ and $G^{p^{\ell}}\in K^{\flat \circ}[x_1,\dots,x_N].$ By Lemma \ref{lemgoodrepapro}, we may suppose that $G_{c}^{p^{\ell}}=\sum_{i\geq 0}^{m}u^ig_{c,i}$ where $g_{c,j}\in \overline{\F_p}[x_1,\dots,x_N]$, $u\in \overline{\F_p((t))}^{\circ}$ with norm $|u|=|t|^{1/[\overline{\F_p}((t))(u):\overline{\F_p}((t))]}$ and $|u|^m> |t|^{(c+1/2)p^{\ell}}$.

Denote by $I_c$ the ideal of $K^{\flat}[x_1,\dots,x_N]$ generated by all $g_{c,i}$.

If $x\in R(U_0^{\flat})$ is a point such that for all $g\in I_c$, $g(x)=0$, then we have $$|H(\pi(\rho^{-1}((\pi^{b})^{-1}(x))))|=|H(\pi(\rho^{-1}((\pi^{b})^{-1}(x))))-G_c^{\#}(\rho^{-1}((\pi^{\flat})^{-1}(x)))|$$
$$\leq  \max\{|p|^{1/2}|G_c^{\#}(\rho^{-1}((\pi^{\flat})^{-1}(x))|,|p|^{c+1/2}\}=|p|^{c+1/2}.$$


On the other hand we have the following lemma.
\begin{lem}\label{lemperinimpeq}Let $x$ be a point in $\Per\cap U_0^{\ad}$ satisfying $|H(x)|\leq |p|^{c+1/2}$. Then for all $g\in I_c$, we have $g(\phi(x))=0.$
\end{lem}
\proof[Proof of Lemma \ref{lemperinimpeq}]Observe that $|(H-G_{c}^{\#})(\chi(x))|\leq |p|^{1/2}\max(|H(x)|,|p|^c)$ and $|H(x)|\leq |p|^{c+1/2}$. We have
$|G_{c}(\phi(x))|=|G_{c}^{\#}(\chi(x))|\leq |p|^{c+1/2}.$ For all $j\geq 0$, we have $g_{c,j}(\phi(x))\in k$. It follows that either $|g_{c,j}(\phi(x))|=1$ or $g_{c,j}(\phi(x))=0$ for all $j\geq 0.$ If $g_{c,j}(\phi(x))=0$ for all $j\geq 0$, then for all $g\in I_c$ we have $g(\phi(x))=0.$
Otherwise, let $j_0$ be the smallest $j$ satisfying $|g_{c,j}(\phi(x))|=1.$
It follows that $|G_c(\phi(x))|=|u|^{j_0/p^{\ell}}.$ Since $|G_{c}(\phi(x))|=|G_{c}^{\#}(\chi(x))|\leq |p|^{c+1/2}$, we get a contradiction.
\endproof

Set $I:=\sum_{c\geq 1}I_c$. Since $K^{\flat}[x_1,\dots,x_N]$ is Noetherian, there exists $M\in \Z^+$, such that $I=\sum_{c=1}^{M}I_c.$ Set $\varepsilon:=|p|^{M+1/2}$.
Let $x$ be a point in $\Per\cap U_0^{\ad}$ satisfying $|H(x)|\leq \varepsilon=|p|^{M+1/2}$.
By Lemma \ref{lemperinimpeq}, for all $g\in I=\sum_{c=1}^{M}I_c$, we have $|g(\phi(x))|=0.$
It follows that for all $c\geq 1$ and $g\in I_c$, we have $|g(\phi(x))|=0$. Then we have $|H(x)|\leq |p|^{c+1/2}$ for all $c\geq 0.$ Let $c$ tend to infinity, we have $H(x)=0$. We conclude our proof.
\endproof

\bigskip

\subsection{Proof of Theorem \ref{thmmmdlfper}}
Suppose that $V\cap \Per$ is Zariski dense in $V$.
We claim the following
\begin{lem}\label{lemsubperzari}There exists a  Zariski dense 
subset $S \subseteq V$ with the property  that $F^{\ell}(S)= S$ for some  positive integer ${\ell}$.
\end{lem}

Since $S$ is Zariski dense in $V$ and $S=F^{{\ell}}(S)$ is Zariski dense in $F^{{\ell}}(V)$. It follows that $V=F^{{\ell}}(V).$
Then Lemma \ref{lemsubperzari} implies Theorem \ref{thmmmdlfper}.

\medskip

Now we prove Lemma \ref{lemsubperzari} in the rest of this section.

Since $\cup_{i=0}^{N}\tau^{-1}(\Per_i^{\ad}\cap V_i^{\ad})=\Per$ is Zariski dense in $V$, there exists $i=0,\dots,N$, such that $\tau^{-1}(\Per_i^{\ad}\cap V_i^{\ad})$ is Zarisi dense in $V.$ We may suppose that $i=0$.

Let $Z$ be the Zariski closure of $\phi(\tau^{-1}(\Per_0^{\ad}\cap V_0^{\ad}))\subseteq \P^N_{K^{\flat}}.$ Since $\phi(\tau^{-1}(\Per_0^{\ad}\cap V_0^{\ad}))$ is defined over $k$ and it is Zariski dense in $Z$, $Z$ is defined over $k=\overline{\F_p}$. Then $Z$ is defined over a finite extension of $\F_p$. It follows that there exists ${\ell}\geq 1$, such that $\Phi^{sl}(Z)=Z$.

Set $S^{\flat,\ad}:=\tau^{\flat}(Z(K^{\flat}))\cap U_0^{\flat}$. We have $\iota(\Per^{\ad}\cap V_0^{\ad})\subseteq S^{\flat,\ad}\cap \pi^{\flat}(\rho(\pi^{-1}(V^{\ad}))).$

We claim the following
\begin{lem}\label{lemsbinvb}We have $S^{\flat,\ad}\subseteq \pi^{\flat}(\rho(\pi^{-1}(V_0^{\ad})))$.
\end{lem}

\begin{rem}
We note that if $\pi^{\flat}(\rho(\pi^{-1}(V_0^{\ad})))$ is algebraic, then our lemma is easy. Since $\phi(\tau^{-1}(\Per_0^{\ad}\cap V_0^{\ad}))$ is Zariski dense in $Z$, and $\pi^{\flat}(\rho(\pi^{-1}(V_0^{\ad})))$ is algebraic, we have $S^{\flat,\ad}\subseteq \pi^{\flat}(\rho(\pi^{-1}(V_0^{\ad})))$.

But in general $\pi^{\flat}(\rho(\pi^{-1}(V_0^{\ad})))$ is not algebraic since the map $\rho$ is not algebraic. Our proof of Lemma \ref{lemsbinvb} is based on Lemma \ref{lemappro}, which allows us to approximate $\pi^{\flat}(\rho(\pi^{-1}(V_0^{\ad})))$ by algebraic subvarieties.
\end{rem}

By assuming Lemma \ref{lemsbinvb}, we have $\pi(\rho^{-1}((\pi^{\flat})^{-1}(S^{\flat,\ad})))\subseteq V_0^{\ad}$. Set $S=\tau^{-1}(\pi(\rho^{-1}((\pi^{\flat})^{-1}(S^{\flat,\ad}))))$. We have $S\subseteq V$ is a Zariski dense subset of $V$. Moreover, we have $F^{{\ell}}(S)=S$. This concludes the proof of Lemma \ref{lemsubperzari}.

\medskip

Now we only need to prove Lemma \ref{lemsbinvb}.

At first, we need the following
\begin{lem}\label{lemconzaricona}Let $H\in K^{b}[z_{0,1},\dots,z_{0,N}]$ be a polynomial with norm $1$. Suppose that for every point $x\in \iota(\Per^{\ad}\cap V_0^{\ad})$, we have $|H(x)|\leq 1/p^s$, where $s\in \Z^+$. Then for every point $y\in S^{\flat,\ad}$, we have $|H(y)|\leq 1/p^s$.
\end{lem}

\proof[Proof of Lemma \ref{lemconzaricona}]Observe that we have a map $$R(U_0^{\ad})= (K^{\flat\circ})^N\to (K^{\flat\circ}/(t^s))^N=\A^N_{K^{\flat\circ}/(t^s)}(K^{\flat\circ}/(t^s))$$
defined by $(x_1,\dots,x_n)\mapsto(\overline{x_1},\dots,\overline{x_N})$ where $\overline{x_i}=x_i \mod t^s$.
Denote by $$\overline{H}:=H \mod t^s.$$ For every point $x\in \iota(\Per^{\ad}\cap V_0^{\ad})$, we have $\overline{H}(\overline{x})=0.$
Observe that
$$\overline{\iota(\Per^{\ad}\cap V_0^{\ad})}=(\phi(\tau^{-1}(\Per_0^{\ad}\cap V_0^{\ad})))\times_{\Spec k}\Spec (K^{\flat,\circ}/(t^s))$$
is Zariski dense in $Z\times_{\Spec k}\Spec (K^{\flat,\circ}/(t^s)).$ It follows that $$Z\times_{\Spec k}\Spec (K^{\flat,\circ}/(t^s))\subseteq \{\overline{H}=0\}.$$
Then we have $$\overline{S^{\flat,\ad}}=Z(k)\times_{\Spec k}\Spec (K^{\flat,\circ}/(t^s))\subseteq Z\times_{\Spec k}\Spec (K^{\flat,\circ}/(t^s))\subseteq \{\overline{H}=0\}.$$
It follows that for every $x\in S^{\flat,\ad}$, we have $H(x)=0 \mod t^s$. Then we have $|H(x)|\leq 1/p^s$, for all $x\in S^{\flat,\ad}$.
\endproof

Now we apply Lemma \ref{lemappro} to $H_{0,j}\in K\langle z_{0,1},\dots,z_{0,N}\rangle\subseteq R_0^{\perf}$ for every $j=1,\dots,m$. For any $s\geq 2$
there exists $h_{s}\in R^{\flat,\perf}_{0}$ such that for all $x\in U_0^{\perf}$, we have
$$|H_{0,j}(x)-h_s^{\#}(x)|\leq |t|^{1/2}\max(|H_{0,j}(x)|,|t|^s)=|t|^{1/2}\max\{|h_s^{\#}(x)|,|t|^s\}<1 \eqno (2).$$
It follows that $\|h_s\|=\|H_{0,j}\|=1$.

For every point $x^{\flat}\in (\pi^{\flat})^{-1}(\iota(\Per^{\ad}\cap V_0^{\ad}))$, we have $$x:=\rho^{-1}(x^{\flat})\in \pi^{-1}(\Per^{\ad}\cap V_0^{\ad}).$$ Then we have $H_{0,j}(x)=0.$ By $(2)$, we have
$$|h_s(x^{\flat})|\leq |t|^{s+1/2}=|t|^{1/2}\max\{|h_s(x^{\flat})|,|t|^s\}=1/p^{s+1/2} .$$
Since $h_s\in R_0^{\flat,\perf}=K^{\flat}\langle z_{0,1}^{1/p^{\infty}},\dots,z_{0,N}^{1/p^{\infty}} \rangle$, there are $r\geq 0$ and a function
$$g_s\in K^{\flat}[ z_{0,1}^{1/p^{r}},\dots,z_{0,N}^{1/p^{r}}]$$ such that $\|h_s-g_s\|< 1/p^{s}.$ It follows that $g_s^{p^r}\in K^{\flat}[ z_{0,1},\dots,z_{0,N}]$ and $$\|h_s^{p^r}-g_s^{p^r}\|\leq |p|^{sp^r}.$$

Then for every point $x^{\flat}\in (\pi^{\flat})^{-1}(\iota(\Per^{\ad}\cap V_0^{\ad}))$, we have
$$|g_s^{p^r}(\pi^{\flat}(x^{\flat}))|=|g_s^{p^r}(x^{\flat})|=|h_s^{p^r}(x^{\flat})+(g_s^{p^r}(x^{\flat})-h_s^{p^r}(x^{\flat}))|\leq |p|^{sp^r}.$$
By Lemma \ref{lemconzaricona}, for all $y\in S^{\flat,\ad}$, we have $|g_s^{p^r}(y)|\leq |p|^{sp^r}$. Then we have
$|h_s((\pi^{\flat})^{-1}(y)|\leq 1/p^{s}$ and
$$|h_s^{\#}(\rho^{-1}((\pi^{\flat})^{-1}(y))|=|h_s((\pi^{\flat})^{-1}(y)|\leq 1/p^{s}$$ for all $y\in S^{\flat,\ad}$.

By Equation $(2)$, we have
$$|H_{0,j}(x)-h_s^{\#}(x)|\leq |t|^{1/2}\max\{|h_s^{\#}(x)|,|t|^s\}= 1/p^{s+1/2}$$ for all $x\in \rho^{-1}((\pi^{\flat})^{-1}(S^{\flat,\ad})).$
It follows that for all $x\in \rho^{-1}((\pi^{\flat})^{-1}(S^{\flat,\ad})),$ we have $|H_{0,j}(x)|\leq 1/p^s$. Let $s\to \infty$, we have $|H_{0,j}(x)|=0$ for all $x\in \rho^{-1}((\pi^{\flat})^{-1}(S^{\flat,\ad})).$ Since $|H_{0,j}(x)|=|H_{0,j}(\pi(x))|$, we have $|H_{0,j}(y)|=0$ for all $j=1,\dots,m$ and $y\in \pi(\rho^{-1}((\pi^{\flat})^{-1}(S^{\flat,\ad}))).$ It follows that $\pi(\rho^{-1}((\pi^{\flat})^{-1}(S^{\flat,\ad})))\subseteq V_0^{\ad}.$ Then we have $S^{\flat,\ad}\subseteq \pi^{\flat}(\rho(\pi^{-1}(V_0^{\ad})))$.
\endproof
\subsection{Scanlon's proof of Theorem \ref{thmmmdlfper}} In this section, we discuss Scanlon's proof of Theorem \ref{thmmmdlfper}. In this proof, we don't need the perfectoid spaces.

\medskip

Let $V$ be a subvariety of $\P^N$ such that $\Per \cap V$ is Zariski dense in $V$. We want to show that $V$ is periodic.

\medskip

We first treat the case where $F$ is defined over $\overline{\Q_p}^{\circ}$. Since all points in $\Per$ are defined over $\overline{\Q_p}$ and $\Per \cap V$ is Zariski dense in $V$, $V$ is defined over $\overline{\Q_p}$.
There exists a finite extension $K_p$ of $\Q_p$ such that $F$ is defined over $K_p$ i.e.
$F$ takes form $$F:[x_0:\dots:x_N]\mapsto [x_0^q+p'P_0(x_0,\dots,x_n):\dots: x_N^q+p'P_N(x_0,\dots,x_N)]$$ where $p'\in K_p^{\circ\circ}$, $q$ is a power of $p$, $P_0,\dots,P_N$ are homogeneous polynomials of degree $q=p^s$ in $K_p^{\circ}[x_0,\dots,x_N].$  After replacing $F$ by a suitable iterate, we may assume that the residue field $\tilde{K}:=K^{\circ}/K^{\circ\circ}$ is fixed by the $q$-power Frobenius.

\medskip

By the structure of the absolute Galois group of $K_p$, there exists an element $\sigma\in \Gal(\overline{K_p}/K_p)$ which lifts the $q$-power Frobenius.
Then we have the following lemma.
\begin{lem}[\cite{Medvdev}]\label{lemsceqdi}We have $\Per=\{x\in \P^N(\overline{\Q_p})|\,\, F(x)=\sigma(x)\}.$
\end{lem}

\proof[Proof of Lemma \ref{lemsceqdi}]
Recall that the reduction map $${\red}:\P^N(\overline{\Q_p})\to \P^N(\overline{\F_p})$$ gives a bijection between $\Per$ and $\P^N(\overline{\F_p}).$

Let $x$ be any point in $\Per$. We have that $F(x)\in \Per$ and $\red(F(x))=\red(x)^q$. On the other hand, we have that $\sigma(x)\in \Per$ and $\red(\sigma(x))=\red(x)^q.$ Then we have $F(x)=\sigma(x).$

Let $x$ be any point in $\P^N(\overline{\Q_p})$ satisfying $F(x)=\sigma(x)$. Since $x$ is defined over a finite extension of $K_p$, there exists $n\geq 1$ such that $\sigma^n(x)=x$. It follows that $$F^n(x)=F^{n-1}(\sigma(x))=\sigma(F^{n-1}(x))=\dots=\sigma^n(x)=x.$$
Then $x$ is periodic.
\endproof

Observe that $\sigma(V)$ is a subvariety of $\P^N.$ Then we have $$\sigma(V\cap \Per)=F(V\cap \Per)\subseteq \sigma(V)\cap F(V).$$

Since $V\cap \Per$ is Zariski dense in $V$, we have $\sigma(V)=F(V).$ Since $V$ is defined over a finite extension of $\Q_p$, there exists $n\geq 1$ such that $\sigma^n(V)=V.$ It follows that $$F^n(V)=F^{n-1}(\sigma(V))=\sigma(F^{n-1}(V))=\dots=\sigma^n(V)=V.$$ Then $V$ is periodic.

\medskip

Now we treat the general case.

There exists subring $R\subseteq \C_p^{\circ}$ which is finitely generated and $\Z$ such that $F$ is defined over $R$.
Let $m:=R\cap \C_p^{\circ\circ}$ be a maximal ideal of $R.$
By Lemma \ref{lemembedalg},  there exists $\sigma\in \Gal(\C_p/\Q)$ such that
$\sigma(R)\subseteq \overline{\Q_p}^{\circ}\subseteq \C_p^{\circ}$ and $\sigma(m)=\overline{\Q_p}^{\circ\circ}\cap R.$

Denote by $F^{\sigma}$ the Galois conjugate of $F$ by $\sigma$ i.e. $F^{\sigma}$ is obtained by changing every coefficient of $F$ by its image under $\sigma.$
Since $F^{\sigma} \mod \C_p^{\circ\circ}=F \mod \C_p^{\circ\circ}$, $F^{\sigma}$ is a lift of Frobenius on $\P^N_{\C_p}.$ Moreover it is defined over $\overline{\Q_p}^{\circ}$. 

Since $V\cap \Per$ is Zariski dense in $V$,  $\sigma(V)\cap \sigma(\Per)$ is Zariski dense in $\sigma(V)$. 
Moreover $\sigma(\Per)$ is exactly the set of periodic points of $F^{\sigma}.$ 
Then the previous argument shows that $\sigma(V)$ is Periodic under $F^{\sigma}.$ It follows that $V$ is periodic under $F$.

\bigskip

\section{Coherent backward orbits}\label{sectioninverseorbit}
In this section, we let $K=\C_p$. Then $K$ is a perfectoid field and $K^{\flat}$ is the completion of the algebraical closure of $\F_p((t)).$ We may suppose that $|p|=|t|=p^{-1}.$ Let $k=\overline{\F_p}$ which is a subfield of $K^{\flat}.$

Let $F:\P^N_{K}\rightarrow \P^N_{K}$ be an endomorphism taking form $$F:[x_0:\dots:x_N]\mapsto [x_0^q+p'P_0(x_0,\dots,x_N):\dots: x_N^q+p'P_N(x_0,\dots,x_N)]$$ where $p'\in K^{\circ\circ}$, $q$ is a power of $p$, and $P_0,\dots,P_N$ are homogeneous polynomials of degree $q$ in $K^{\circ}[x_0,\dots,x_N].$

\medskip

The aim of this section is to prove Theorem \ref{thmrevedml} and Theorem \ref{thmrevedmllimit}.

\medskip


Without loss of generality, we may suppose that $b_0\in R(U_0^{\ad})$. It follows that $b_i\in R(U_0^{\ad})$ for all $i\geq 0.$
Set $w:=\pi^{\flat}\circ\rho\circ\psi^{-1}((b_0,b_1\dots))\in \P^N_{K^{\flat}}(K_b).$ Then $w\in R(U_0^{\flat,\ad}):=\{[1:x_1:\dots:x_N]|\,\, |x_i|\leq 1\}\subseteq \P^N_{K^{\flat}}(K_b).$ It follows that $w^{1/q^{n}}\subseteq R(U_0^{\flat,\ad})$ for all $n\geq 0.$

If $\{b_i\}_{i\geq 0}$ is infinite, we may suppose that $b_1\neq b_0$ and then $b_i$, $i\geq 0$ are all different.
Let $Z$ be the reduced subvariety of  $U_0^{\flat}:=\Spec K^{\flat}[x_1,\dots,x_N]$,  whose support is the union of all positive dimensional irreducible components of the Zariski closure of $\{w^{1/q^{n_i}}\}_{i\geq 0}.$

There exists $A\geq 0$, such that $Z$ is the Zariski closure of $\{w^{1/q^{n_i}}\}_{i\geq A}$ in $U_0^{\flat}$. Moreover, for all $n\geq A$, $Z$ is the  Zariski closure of $\{w^{1/q^{n_i}}\}_{i\geq n}$ in $U_0^{\flat}$.

Denote by $I(Z)$ the ideal in $K^{\flat}[x_1,\dots,x_N]$ which defines $Z$.

For every polynomial $f=\sum_{I}a_Ix^I\in K^{\flat}[x_1,\dots.x_N]$ and $i\in \Z$, we denote by $f^{\sigma^i}:=\sum_{I}a_I^{q^i}x^I.$ Observe that $f(y^{1/{q^i}})=(f^{\sigma^i}(y))^{1/q^i}$ for all $i\geq 0$ and $y\in R(U_0^{\flat,\ad}).$

Then we have the following lemma
\begin{lem}\label{lemneartheneq}Let $f\in k[x_1,\dots,x_N]$ be a polynomial defined over $k$. If there exists $c\in (0,1)$ and $B\geq A$, such that for all $i\geq B$, $|f(w^{1/q^{n_i}})|\leq c$, then $f\in I(Z).$
\end{lem}
\proof[Proof of Lemma \ref{lemneartheneq}]
There exists $L\geq 1$ such that $f$ is defined over $\F_{q^L}.$ Then we have $f^{\sigma^{nL}}=f$ for all $n\geq 0.$
For $t=0,\dots,L-1$, set $T_t:=\{i\geq B|\,\, n_i=t \mod L\}.$

For all $t=0,\dots,L-1$ satisfying $\#T_t=\infty$, we have $$|f(w^{1/q^{t}})|^{1/q^{n_i-t}}=|f^{\sigma^{n_i-t}}(w^{1/q^{t}})|^{1/q^{n_i-t}}=|f(w^{1/q^{n_i}})|\leq c,$$ for all $i\in T_t.$ It follows that $|f(w^{1/q^{t}})|\leq c^{q^{n_i-t}}$ for all $i\in T_t.$ Since $T_t$ is infinite, $n_i$ can be arbitrary large. Then we have  $|f(w^{1/q^{t}})|=0$ for all $i\in T_t.$ It follows that $$|f(w^{1/q^{n_i}})|=|f^{\sigma^{n_i-t}}(w^{1/q^{t}})|^{1/q^{n_i-t}}=|f(w^{1/q^{t}})|^{1/q^{n_i-t}}=0$$
for all $i\in T_t.$ 
Set $$T':=\sqcup_{0\leq t\leq L-1,\#T_t=\infty}T_t.$$
It follows that $f(w^{1/q^{n_i}})=0$ for all $i\in T'.$ Since $\{i\geq A\}\setminus T'$ is finite, $\{w^{1/q^{n_i}}\}_{i\in T'}$ is Zariski dense in $Z$. Then $f\in I(Z).$
\endproof

\begin{lem}\label{lemzoverfini}We have that $Z$ is defined over $k$. In particular, there exists $r\geq 1$ such that $\Phi^{sr}(Z)=Z$ and $\{w^{1/q^{i}}\}_{i\in \Z}\subseteq \cup_{i=0}^{r-1}\Phi^{si}(Z).$
\end{lem}

\proof[Proof of Lemma \ref{lemzoverfini}]We only need to show that $I(Z)$ is generated by finitely many polynomials in $k[x_1,\dots,x_N]\subseteq K^{\flat}[x_1,\dots,x_N]$.
In fact, if $I(Z)=(g_1,\dots,g_l)$ and $g_i\in k[x_1,\dots,x_N]$ for all $i=1,\dots,{\ell},$ then there exists $r\geq 1$ such that all the coefficients of  $g_i, i=1,\dots,{\ell}$ are defined over $\F_{q^r}.$ Then we have $\Phi^{sr}(Z)=Z.$ Moreover, there exists $j\geq 0$, such that $w^{1/p^j}\in Z$. It follows that $\{w^{1/q^{i}}\}_{i\in \Z}\subseteq \cup_{i\in\Z}\Phi^{si}(Z)=\cup_{i=0}^{r-1}\Phi^{si}(Z).$

Write $I(Z)=(f_1,\dots, f_m)$ where $m\geq 1$ and $f_i\in K^{\flat}[x_1,\dots,x_N]$ for all $i=1,\dots,m.$
Denote by $d:=\max_{0\leq i\leq m}\{\deg(f_i)\}.$

By Lemma \ref{lemgoodrepapro}, for all $i=1,\dots,m$, there exists a sequence of polynomial $\{f_{i,n}\}_{n\geq 1}$ such that $\|f_i-f_{i,n}\|\leq |t^n|$ and taking form
$f_{i,n}=\sum_{j=0}^{m_{i,n}}u_{i,n}^jf_{i,n,j}$ where $f_{i,n,j}\in \overline{\F_p}[x_1,\dots,x_N]$ of degree at most $d$, $u_{i,n}\in \overline{\F_p((t))}^{\circ}$ with norm $|u_{i,n}|=|t|^{1/[\overline{\F_p}((t))(u_{i,n}):\overline{\F_p}((t))]}$ and $|u_{i,n}|^{m_{i,n}}> |t|^{n}$.

\medskip

We claim that $f_{i,j,n}\in I(Z)$ for all $j=0,\dots,m_{i,n}.$

We prove that claim by induction on $j.$ For $j=0$, we have
$$|f_{i,0,n}(w^{1/q^{n_l}})|=|f_{i,n}(w^{1/q^{n_l}})-\sum_{j\geq 1}u_{i,n}^{j}f_{i,j,n}(w^{1/q^{n_l}})|$$
$$\leq \max\{|t^n|,|u_{i,n}|\}<1$$ for all ${\ell}\geq A.$
By Lemma \ref{lemzoverfini}, we have $f_{i,0}\in I(Z).$

If $j\geq 1$ and $f_{i,0,n},\dots,f_{i,j-1,n}\in I(Z)$, then
$$|f_{i,j,n}(w^{1/q^{n_l}})|=|u_{i,n}^{-j}(f_{i,n}(w^{1/q^{n_l}})-\sum_{0\leq t'\leq j-1}u_{i,n}^{t'}f_{i,t',n}(w^{1/q^{n_l}}))-\sum_{t'\geq j+1}u_{i,n}^{t'-j}f_{i,j,n}(w^{1/q^{n_l}})|$$
$$\leq \max\{|t|^n/|u_{i,n}|^{j}, |\sum_{t'\geq j+1}u_{i,n}^{t'-j}f_{i,j,n}(w^{1/q^{n_l}})|\}\leq \max\{|t|^n/|u_{i,n}|^{j}, |u_{i,n}|\}<1$$ for all ${\ell}\geq A.$ By Lemma \ref{lemzoverfini}, we have $f_{i,j,n}\in I(Z).$
Then we conclude the proof of the claim. It follows that $f_{i,n}\in I.$

Set $I_d:=\{f\in I|\,\, \deg(f)\leq d\}$. Then $I_d$ is a finite-dimensional $K^{\flat}$-vector space.
For all $n\geq 0, j=0,\dots,m_{i,n}$, denote by $I_{i,j,n}$ the $K_b$-vector space spanned by $f_{i,0,0}\dots,f_{i,j,0},\dots, f_{i,0,n}\dots,f_{i,j,n}.$ Then $\cup_{n\geq 0, j=0,\dots,m_{i,n}}I_{i,j,n}$ is a subspace of $I_d.$ Since $\dim I_d$ is finite, $\cup_{n\geq 0, j=0,\dots,m_{i,n}}I_{i,j,n}$ is closed. Observe that $f_i$ is contained in the closure of $\cup_{n\geq 0, j=0,\dots,m_{i,n}}I_{i,j,n}$, we have $f_i\in \cup_{n\geq 0, j=0,\dots,m_{i,n}}I_{i,j,n}.$ There exists $l_i\geq 0$, such that $f_i\in I_{i,m_{i,l_i},l_i}.$ It follows that $I=(f_1,\dots,f_m)\subseteq \sum_{1\leq i\leq m}(I_{i,m_{i,l_i},l_i})\subseteq I$. Then we have $I=(f_{i,j,n})_{1\leq i\leq m, 0\leq n\leq l_i,0\leq j\leq m_{l_i}}$ and $f_{i,j}\in k[x_1,\dots,x_N]$ for all $i,j.$
\endproof

\subsection{Proof of Theorem \ref{thmrevedmllimit}}
Let $V$ be a subvariety of $\P^N_{\C_p}$ such that there exists a subsequence $\{b_{n_i}\}_{i\geq 0}$ such that $|d(b_{n_i},V)|\to 0$ when $i\to\infty$. We need to show that $b_{n_i}\in V$ for $i$ large enough and there exists $r\geq 0$, such that $\{b_i\}_{i\geq 0}\subseteq \cup_{i=0}^{r-1}F^i(V).$

If $\{b_i\}_{i\geq 0}$ is finite, Theorem \ref{thmrevedmllimit} is trivial. So we suppose that $\{b_i\}_{i\geq 0}$ is infinite.

Denote by $I(V)$ the ideal in $K[x_1,\dots,x_N]$ which defines $V\cap U_1.$
Then for any point in $R(U^{\ad}_0)$, we have $d(y,V)=\max\{|H(y)||\,\, H\in I(V) \text{ and } \|H\|=1\}.$

\bigskip

Denote by $Z$ the union of all positive dimensional irreducible components of the Zariski closure of $\{w^{1/q^{n_i}}\}_{i\geq 0}.$ There exists $A\geq 0$, such that $Z$ is that Zariski closure of $\{w^{1/q^{n_i}}\}_{i\geq A}$ in $U_0^{\flat}:=\Spec K^{\flat}[x_1,\dots,x_N]$. Moreover, for all $n\geq A$, $Z$ is the Zariski closure of $\{w^{1/q^{n_i}}\}_{i\geq n}$ in $U_0^{\flat}$.

Denote by $I(Z)$ the ideal in $K^{\flat}[x_1,\dots,x_N]$ which defines $Z$.

Let $H$ be a polynomial in $I(V).$

\begin{lem}\label{lemzinhzero1}For any point $x\in Z\cap R(U_0^{\flat,\ad}),$ we have $H(\pi(\rho^{-1}((\pi^{\flat})^{-1}(x))))=0$.
\end{lem}

\proof[Proof of Lemma \ref{lemzinhzero1}]By Remark \ref{remappro}, for any $c \in \Z^+$, there exists ${\ell}\in \N$ and an element $G_c\in K^{\flat \circ}[x_1^{1/p^{{\ell}}},\dots, x_N^{1/p^{{\ell}}}]$ such that for all $x\in U_0^{\perf}$, we have
$$|H\circ \pi(x)-G_{c}^{\#}(x)|\leq |p|^{1/2}\max(|H(x)|,|p|^c)=|p|^{1/2}\max(|G_{c}^{\#}(x)|,|p|^c),$$ and $G_c^{p^{\ell}}\in K^{\flat \circ}[x_1,\dots,x_N].$
By Lemma \ref{lemgoodrepapro}, we may suppose that $G_{c}^{p^{\ell}}=\sum_{i\geq 0}^{m}u^ig_{i}$ where $g_{i}\in \overline{\F_p}[x_1,\dots,x_N]$, $u\in \overline{\F_p((t))}^{\circ}$ with norm $|u|=|t|^{1/[\overline{\F_p}((t))(u):\overline{\F_p}((t))]}$ and $|u|^m> |t|^{(c+1/2)p^{\ell}}$.

There exists $A_1\geq 0$, such that for all $i\geq A_1$, $|H(b_{n_i})|\leq |p|^{c+1}.$

For all $i\geq A_1$, we have $$|G_{c}(w^{1/q^{n_i}})|\leq \max\{|H(b_{n_i})|,|H(b_{n_i})-G_{c}^{\#}(\rho^{-1}(w^{1/q^{n_i}}))|\}\leq |p|^{c+1/2}.$$
Then we have $$|G_{c}(w^{1/q^{n_i}})^{p^{\ell}}|\leq |t|^{p^{\ell}(c+1/2)}.$$

We claim that for all $j=0,\dots,m$, we have $g_j\in I(Z).$

We prove this claim by induction on $j$. Suppose that for all $0\leq t'<j\leq m$, we have $g_{t'}\in I(Z).$ For all $i\geq \max\{A, A_1\}$, we have
$$ |u^jg_j(w^{1/q^{n_i}})+\sum_{t'\geq j+1}u^{t'}g_{t'}(w^{1/q^{n_i}})|=|G_{c}(w^{1/q^{n_i}})^{p^{\ell}}|\leq |t|^{p^{\ell}(c+1/2)}.$$
It follows that $|g_j(w^{1/q^{n_i}})|\leq \max\{|t|^{p^{\ell}(c+1/2)}/|u|^j, |u|\}<1$ for all $i\geq \max\{A, A_1\}.$ Then Lemma \ref{lemneartheneq} implies that $g_j\in I(Z)$ for $j=0,\dots,m.$ Then we conclude our claim.

Then for any $x\in Z\cap R(U_0^{\flat,\ad}),$ we have
$$|H(\pi(\rho^{-1}((\pi^{b})^{-1}(x))))|=|H(\pi(\rho^{-1}((\pi^{b})^{-1}(x))))-G_c^{\#}(\rho^{-1}((\pi^{\flat})^{-1}(x)))|.$$
$$\leq \max\{|p|^{1/2}|G_c^{\#}(\rho^{-1}((\pi^{\flat})^{-1}(x))|,|p|^{c+1/2}\}=|p|^{c+1/2}$$
%
%
Let $c$ tends to infinity, then we have $|H(\pi(\rho^{-1}((\pi^{b})^{-1}(x))))|=0.$ We complete the proof of our lemma.
\endproof
This lemma shows that $S:=\pi(\rho^{-1}((\pi^{b})^{-1}(Z\cap R(U_0^{\flat,\ad}))))\subseteq V.$ Then $b_{n_i}\in V$ for $i\geq A$. By Lemma \ref{lemzoverfini}, there exists $r\geq 1$ such that $\Phi^{sr}(Z)=Z$ and $\{w^{1/p^{i}}\}_{i\in \Z}\subseteq \cup_{i=0}^r\Phi^{si}(Z).$ It follows that $\{b_i\}_{i\geq 0}\subseteq \cup_{i=0}^{r-1}F^i(S)\subseteq \cup_{i=0}^{r-1}F^i(V).$
\subsection{Proof of Corollary \ref{cordmlrevtv}}Let $V$ be a subvariety of $\P^N_{\C_p}$ of positive dimension. We need to show there exists $c>0$ such that for all $i\geq 0$, either $b_i\in V$ or $d(b_i,V)>c.$

Otherwise, there exists a subsequence $\{b_{n_i}\}_{i\geq 0}\subseteq \{b_i\}_{i\geq 0}\setminus V$ such that $d(b_{n_i},V)$ tends to $0.$ By Theorem \ref{thmrevedmllimit}, we have $b_{n_i}\in V$ for sufficiently large $i$, which is a contradiction.

\subsection{Proof of Theorem \ref{thmrevedml}}
Let $V$ be a positive subvariety of $\P^N_{\C_p}$ such that $\{b_i\}_{i\geq 0}\cap V$ is Zariski dense in $V$. Let $\{n_1<n_2<\dots\}$ be the set of $n\geq 0$ such that $b_n\in V.$
We need to show that $V$ is periodic under $F$.

\bigskip

%
%
%
%
If $\{b_i\}_{i\geq 0}$ is finite, then all points in $\{b_i\}_{i\geq 0}$ are periodic.  Moreover $V$ is a union of finitely many periodic points. Then $V$ is periodic.

Now we may suppose that $\{b_i\}_{i\geq 0}$ is infinte.

Denote by $I(V)$ the ideal in $K[x_1,\dots,x_N]$ which defines $V\cap U_1.$
Let $H$ be a polynomial $I(V).$

By Lemma \ref{lemzinhzero1}, for any point $x\in Z\cap R(U_0^{\flat,\ad}),$ we have $H(\pi(\rho^{-1}((\pi^{b})^{-1}(x))))=0$.

It follows that that $S:=\pi(\rho^{-1}((\pi^{b})^{-1}(Z\cap R(U_0^{\flat,\ad}))))\subseteq V.$ Since $b_{n_i}\in S$ for all $i\geq A,$  $S$ is Zariski dense in $V.$
Since $\Phi^{rs}(Z\cap R(U_0^{\flat,\ad}))=Z\cap R(U_0^{\flat,\ad})$, we have $F^r(S)=S$. It follows that $F^r(V)=V$. Then we conclude the proof.
\endproof

\section{Appendix}

Let $X$ be any projective variety over $\C_p$ and $F:X\to X$ be an endomorphism. 
Let $\mathfrak{X}\to\Spec \C_p^{\circ}$ be a 
finitely presented projective scheme which is flat over $\Spec \C_p^{\circ}$ whose generic fiber is $X$ and $\sL$ an ample line bundle on $\fX$. If there exists an endomorphism $\widetilde{F}$ of $\mathfrak{X}$ over $\C_p^{\circ}$ such that $\wF^*\sL=\sL^{\otimes q}$ where $q=p^s, s\geq 1$, the restriction of $\widetilde{F}$ on the generic fiber is $F$ and the restriction $\bar{F}$ of $\widetilde{F}$ on the special fiber $\bar{X}$ is a power of the Frobenius, then we say that \emph{$F$ is a polarized lift of Frobenius on $X$ w.r.t  $(\fX, \widetilde{F}, \sL).$} In particular,  a lift of Frobenius on $\P^N_{\C_p}$ in the previous sections is a lift of Frobenius on $X$ w.r.t a pair $(\P^N_{\C^{\circ}_p}, \widetilde{F}, O_{\P^N_{\C^{\circ}_p}}(1)).$


\medskip

Now assume that $F$ is a polarized lift of Frobenius on $X$ w.r.t the pair $(\fX, \widetilde{F},\sL)$ and we identify $X$ with the generic fiber of $\fX.$

\smallskip

In this appendix, we show that under a technical condition, the dynamical system $(X, F)$ can be embedded in a lift of Frobenius on $\P^N_{\C_p}$ ( w.r.t some $(\P^N_{\C^{\circ}_p}, \widetilde{F}, O_{\P^N_{\C^{\circ}_p}}(1)).$)

\begin{thm}\label{thmembedfrob}Assume that $\fX$ and $\widetilde{F}$ are defined over $\overline{\Q}_p^{\circ}\subseteq \C_p^{\circ}$. Then there exists $N\geq 1$, a lift of Frobenius $G$ on $\P^N_{\C_p}$ 
and an embedding $\tau: X\hookrightarrow \P^N_{\C_p}$ such that $\tau\circ F^l=G\circ \tau$ for some $l\geq 1.$
\end{thm}
This theorem can be viewed as a version of \cite[Proposition 2.1]{fa} for the lifts of Frobenius.

As an application, it implies the dynamical Manin-Mumford Conjecture and Conjecture \ref{condmlrev}, for any polarized lift of Frobenius on $X$ w.r.t some $(\fX, \widetilde{F}, \sL).$

\begin{cor}\label{cormmxcaibliftfro} Let $V$ be any positively dimensional irreducible subvariety of $X$. Denote by $\Per_F$ the set of periodic closed points in $X$. 
Let $\{b_i\}_{i\geq 0}$ be a sequence of closed points in $X$ satisfying $f(b_i)=b_{i-1}$ for all $i\geq 1$.  Then we have that
\begin{points}
\item  if $V\cap\Per_F$ is Zariski dense in $V$, then $V$ is periodic;
\item  if the $\{b_i\}_{i\geq 0}\cap V$ is Zariski dense in $V$, then $V$ is periodic under $F$.
\end{points}
\end{cor}

\proof[Proof of Corollary \ref{cormmxcaibliftfro}]
There exists subring $R\subseteq \C_p^{\circ}$ which is finitely generated over $\Z$ such that $\fX$, $\wF$ and $\sL$ are defined over $R$ i.e. 
there exists a projective scheme $\fX_R$ over $\Spec R$ with an endomorphism $\wF_R$ and an ample line bundle $\sL_R$ such that $\fX=\fX_R\otimes_R\C_p^{\circ}$, $\sL=\sL_R\otimes_R\C_p^{\circ}$, $\wF=\wF_R\otimes_R \C_p^{\circ}$ and $\wF^*_R\sL_R=\sL_R^{\otimes q}.$

If $R\subseteq \overline{\Q_p}^{\circ}$, Theorem \ref{thmembedfrob} reduces it to the case where $X=\P^N_{\C_p}$ and $F$ is a lift of Frobenius on $\P^N_{\C_p}.$
Then we conclude the proof by Theorem \ref{thmmmdlfper} and Theorem \ref{thmrevedmllimit}.

Now we assume that  $R\not\subseteq \overline{\Q_p}^{\circ}$. Set $m:=R\cap \C_p^{\circ\circ}.$ It is a maximal ideal of $R.$
\begin{lem}\label{lemembedalg}Let $R$ be a subring of $\C_p^{\circ}$  which is finitely generated over $\Z$. Let $m:=R\cap \C_p^{\circ\circ}$ be a maximal ideal of $R.$ Then there exists $\sigma\in \Gal(\C_p/\Q)$ such that
$\sigma(R)\subseteq \overline{\Q_p}^{\circ}\subseteq \C_p^{\circ}$ and $\sigma(m)=\overline{\Q_p}^{\circ\circ}\cap R.$
\end{lem}

Now consider $\fX^{\sigma}:=\fX_R\otimes_R^{\sigma}\C_p^{\circ}$, $\sL^{\sigma}:=\sL_R\otimes_R^{\sigma}\C_p^{\circ}$,
$\wF^{\sigma}:=\wF\otimes_R^{\sigma}\C_p^{\circ}$, $X^{\sigma}:=X_R\otimes_R^{\sigma}\C_p$ and $F^{\sigma}:=X_R\otimes_R^{\sigma}\C_p.$
In the tensor product $\d \otimes_R^{\sigma}\C_p$, we use the embedding $\sigma|_R$.
We note that if we view $\C_p$ as an abstract field, $(X, F)$ and $(X^{\sigma}, F^{\sigma})$ are Galois conjugate.
Since the statements of $(i),(ii)$ are purely algebraic, we only need to show it for $(X^{\sigma}, F^{\sigma}).$
Observe that the special fiber $\bar{X}^{\sigma}$ of $\fX^{\sigma}$ is $$\bar{X}^{\sigma}=\fX_R\otimes_R^{\sigma}(\C_p^{\circ}/\C_p^{\circ\circ})=\fX_R\otimes_R(R/m)\otimes^{\bar{\sigma}}_{R/m} (\C_p^{\circ}/\C_p^{\circ\circ})=\fX_R\otimes_R^{\bar{\sigma}}(\C_p^{\circ}/\C_p^{\circ\circ})\simeq\bar{X}.$$
Moreover the restriction of $F^{\sigma}$ on $\bar{X}^{\sigma}$ is exactly $\bar{F}$ under this identification. So $F^{\sigma}$ is some power of Frobenius and $F$ is a lift of the Frobenius w.r.t. $(\fX^{\sigma}, \wF^{\sigma},\sL^{\sigma}).$ 
Since $(\fX^{\sigma}, \wF^{\sigma},\sL^{\sigma})$ is defined over $\sigma(R)\subseteq \overline{\Q_p}^{\circ},$ 
Theorem \ref{thmembedfrob} reduces it to the case where $X=\P^N_{\C_p}$ and $F$ is a lift of Frobenius on $\P^N_{\C_p}.$
Then we conclude the proof by Theorem \ref{thmmmdlfper} and Theorem \ref{thmrevedmllimit}.
\endproof

\proof[Proof of Lemma \ref{lemembedalg}] Since $\C_p$ is algebraically closed, any embedding $R\hookrightarrow \C_p$ extends to an automorphism in $\Gal(\C_p/\Q).$
We only need to find an embedding $\sigma: R\hookrightarrow \C_p^{\circ}\subseteq \C_p$ satisfying $\sigma(m)\subseteq \C^{\circ\circ}_p.$
Indeed since $\sigma^{-1}(\C^{\circ\circ}_p\cap \sigma(R))$ is a maximal ideal of $R$ which contains $m$, we have $\sigma(m)=\sigma(R)\cap \C_p^{\circ\circ}.$

\medskip

Let $t_1,\dots,t_l\in R$ be a set of generators of $R$ over $\Z$. Let $u_1,\dots,u_s$ be a set of generators of $m.$
Set $Y:=\Spec R$ and $Y_{\C_p}:=\Spec R\otimes_{\Z}\C_p.$ We endow $Y_{\C_p}(\C_p)$ the $p$-adic topology induced by the topology on $\C_p.$
An element $f\in R$ can be viewed as an analytic function on $Y_{\C_p}(\C_p)$.

Denote by $i:R\hookrightarrow \C_p$ the inclusion. It defines a point $o\in Y_{\C_p}(\C_p).$ Set $U:=\{x\in Y_{\C_p}(\C_p)|\,\, |t_i|\leq 1, i=1,\dots,l \text{ and } |u_i|<1, i=1,\dots,s\}.$  Then $U$ is an open neighbourhood of $o.$

For any nonzero $P$ element of $R$, denote by $V_P$ the subscheme of $Y_{\C_p}$ defined by $\{P=0\}.$ Since the set of  nonzero prime ideals is countable, and $Y_{\C_p}(\C_p)$ has a complete metric, $Y_{\C_p}(\C_p)\setminus (\cup_{R\setminus \{0\}} V_P)$ is dense in $Y_{\C_p}(\C_p).$ 
Then there exists a point $y\in U\setminus (\cup_{R\setminus \{0\}} V_P).$ It defines a morphism $\sigma:R\to \C_p$ by $f\mapsto f(y).$ Because $y\in U$, we have $\sigma(t_i)\in \C^{\circ}, i=1,\dots,l$ and $\sigma(u_i)\in \C^{\circ\circ}_p, i=1,\dots,s.$ It follows that $\sigma(R)\subseteq \C^{\circ}_p$ and $\sigma(m)\subseteq \C^{\circ\circ}_p.$ Since $y\not\in  (\cup_{R\setminus \{0\}} V_P)$, $\sigma:R\to \C_p$ is an embedding. This concludes the proof.
\endproof

\subsection{Proof of Theorem \ref{thmembedfrob}}
In this section, we assume that $\fX$, $\widetilde{F}$ and $\sL$ are defined over $\overline{\Q}_p^{\circ}\subseteq \C_p^{\circ}$. 
Since $\fX$ is finitely presented, there exists a finite extension $K$ of $\Q_p$ such that $\fX$, $\wF$ and $\sL$ are defined over $K^{\circ}.$
We note that $R:=K^{\circ}$ is a discrete valuation ring. Set $m:=K^{\circ\circ}$ the maximal ideal of $R$ and $\pi$ a generator of $m.$

There exists a flat and geometrically irreducible projective scheme $\fX_R$ over $\Spec R$ an ample line bundle $\sL_R$ and an endomorphism $\wF_R$ such that $\fX=\fX_R\otimes_R\C_p^{\circ}$ $\sL=\sL_R\otimes_R\C_p^{\circ}$  and $\wF=\wF_R\otimes_R \C_p^{\circ}.$ 
We may assume that $\wF_R^*\sL_R=\sL_R^{\otimes q}.$
We denote by $X_s$ the special fiber of $\fX_R$ and $F_s$ the restriction of $\wF_R$ on $X_s.$ Let $X_K$ be the generic fiber of $\fX_R$ and $F_K$ the restriction of $\wF_R$ on $X_K.$  
Write $L_K:=\sL_R|_{X_K}$ and $L_s:=\sL_R|_{X_s}.$  Since $\sL_R$ is ample, after replacing $\sL_R$ by a suitable power, we may assume that $\sL$ is very ample and the morphisms $$\Psi_R:=H^0(X_R,\sL_R)^{\otimes q}\to H^0(X_R,\sL_R^{\otimes q})$$ 
and $$\Psi_s:=H^0(X_s,\sL_s)^{\otimes q}\to H^0(X_s,\sL_R^{\otimes q})$$ are surjective. Moreover, we may assume that $$H^i(X_R, \sL_R)=H^i(X_R, \sL_R^{\otimes q})=0$$ for all $i\geq 1.$ It follows that the natrual morphisms
$$r_1:H^0(X_R,\sL_R)\otimes_R R/m\to H^0(X_s,\sL_s)$$ and $$r_q:H^0(X_R,\sL_R^{\otimes q})\otimes_R R/m\to H^0(X_s,\sL_s^{\otimes q})$$ are isomorphisms.

Since $\wF_R^*\sL_R=\sL_R^{\otimes q},$ 
it induces morphisms
$$\wF_R^*: H^0(X_R,\sL_R)\to H^0(X_R,\sL_R^{\otimes q})$$ and
$$\wF_s^*: H^0(X_s,\sL_s)\to H^0(X_s,\sL_s^{\otimes q}).$$
We have $r_q\circ \wF_R^*=F_s^*\circ r_1$.

Since $R$ is a discrete valuation ring, and $H^0(X_R,\sL_R)$ has no torsions, $H^0(X_R,\sL_R)$  is a free $R$-module. Let $s_0,\dots s_N$ a basis of $H^0(X_R,\sL_R).$
We note that $$r_q(\wF^*(s_i))=r_1(s_i)^q$$ for $i=0,\dots,N.$ It follows that $$\wF^*(s_i)-\Psi_R(s_i^q)\in mH^0(X_R,\sL_R^{\otimes q})=\pi H^0(X_R,\sL_R^{\otimes q})$$ for $i=0,\dots,N.$
In other words, there exists $g_i\in H^0(X_R,\sL_R^{\otimes q})$ such that $$\wF^*(s_i)=s_i^q+\pi g_i, i=1,\dots,N.$$
Since $\Psi_R$ is surjective, there exists $G_i\in R[x_0,\dots,x_N]$ homogenous of degree $q$ such that $g_i=G_i(s_0,\dots,s_N), i=1,\dots,N.$ It follows that 
$$\wF^*(s_i)=s_i^q+\pi G_i(s_0,\dots,s_N), i=1,\dots,N.$$

Let $G_K:\P^N_{K}\to \P^N_{K}$ be the morphism $$[x_0,\dots,x_N]\mapsto [x_0^{q}+\pi G_0(x_0,\dots,x_N):\dots:x_N^{q}+\pi G_N(x_0,\dots,x_N)].$$
Set $G:=G_K\otimes_K\C_p:\P^N_{\C_p}\to \P^N_{\C_p}.$
It is a lift of Frobenius on $\P^N_{\C_p}.$ 
Let $\tau_K: X\to \P^N_{K}$ be the morphism 
$$x\mapsto [s_0(x):\dots: s_N(x)].$$
Since $L_K$ is very ample, $\tau$ is an embedding. We may check that $G_K\circ\tau_K=\tau_K\circ F_K$.
We conclude the proof by setting $\tau:=\tau_K\otimes_K \C_p.$

\newpage
\bibliography{dd}

\end{document}